\documentclass[aap]{imsart}
\usepackage{latexsym,amssymb,amsfonts,amsmath,amsthm,graphicx}
\newtheorem{thm}{Theorem}

\newtheorem{lemma}{Lemma}
\newtheorem{remark}{Remark}
\newtheorem{proposition}{Proposition}

\def\eopt{\hfill$\Box$}
\newcommand{\e}{\varepsilon}
\newcommand{\n}{\noindent}
\newcommand{\vs}{\vspace{10pt}}
\newcommand{\s}{\vspace{2ex}}
\newcommand{\R}{\mathbb{R}}
\newcommand{\PP}{\mathbb{P}}
\newcommand{\E}{\mathbb{E}}

\begin{document}

\begin{frontmatter}

\title{A Spectral Analysis of the Sequence of Firing Phases in Stochastic Integrate-and-Fire Oscillators}
\runtitle{Spectral Analysis of SIFs}

\author{\fnms{Peter} \snm{Baxendale}\thanksref{t1,m1} \ead[label=e1]{baxendal@usc.edu}}
\and
\author{\fnms{John} \snm{Mayberry}\thanksref{t2,m2} \ead[label=e2]{jmayberry@math.cornell.edu}}

\address{Peter Baxendale \\ Department of Mathematics \\
University of Southern California \\
3620 S. Vermont Avenue \\
Los Angeles, CA 90089-2532, U.S.A. \\
\printead{e1}}

\address{John Mayberry \\
 Department of Mathematics \\
 Cornell University \\
 310 Malott Hall \\
 Ithaca, NY 14853-4201, U.S.A. \\
 \printead{e2}}

\affiliation{University of Southern California \thanksmark{m1} and Cornell University \thanksmark{m2}}
\runauthor{Baxendale and Mayberry}

\thankstext{t1}{Supported in part by NSF grant DMS-05-04853.}
\thankstext{t2}{Supported in part by NSF RTG grant DMS-07-39164.}

\begin{abstract}
Integrate and fire oscillators are widely used to model the generation of action potentials in
neurons. In this paper, we discuss small noise
asymptotic results for a class of stochastic integrate and fire oscillators (SIFs) in which the buildup of membrane potential in the neuron is governed by a Gaussian diffusion process. To analyze this model, we study the asymptotic behavior of the spectrum of the firing phase transition operator. We begin by proving strong versions of a law of large numbers and central limit theorem for the first passage-time of the underlying diffusion process across a general time dependent boundary. Using these results, we obtain asymptotic approximations of the transition operator's eigenvalues. We also discuss connections between our results and earlier numerical investigations of SIFs.

\end{abstract}


\end{frontmatter}

\section{Introduction} \label{sec-intro}

The integrate and fire oscillator is widely used to model the
behavior of the membrane potential in a neuron.  Since its
introduction by Lapicque \cite{lap} in 1907 it has been studied by
many authors, see in particular Stein \cite{ste} and Knight
\cite{kni}.  The neurobiological derivation of the model is
described in Tuckwell \cite{tuc1,tuc2} and a recent review of
activity in the field appears in Burkitt \cite{bur1,bur2}.

In this paper we use the following  model for the stochastic
integrate and fire oscillator (SIF).  After starting at some time
$t_0$ the membrane potential $X_t^\e$ evolves according to the
stochastic differential equation
 \begin{equation} \label{X}
  dX_t^\e =  (-\gamma X_t^\e + I(t))dt + \e \, dW_t
  \end{equation}
until it reaches a (time-dependent) threshold level $g(t)$ at time
  \begin{equation} \label{tau}
    \tau_1^\e = \inf\{t \ge t_0 : X_t^\e = g(t)\}.
  \end{equation}
At the hitting time $\tau_1^\e$ the membrane potential discharges,
producing a voltage spike, and resets at a lower value
   \begin{equation}  \label{reset}
     X_{(\tau_1^\e)^+}^\e = h(t).
   \end{equation}
For $t \ge (\tau_1^\e)^+$ the process $X_t^\e$ follows the SDE
(\ref{X}) until the second hitting time $\tau_2^\e = \inf\{t >
\tau_1^\e: X_t^\e = g(t)\}$, and so on, yielding a sequence
$\{\tau_n^\e: n \ge 1\}$ of hitting times.

Here $W_t$ is a standard one-dimensional Wiener process, and $\e
\ge 0$ determines the intensity of the noise in the integrate and
fire oscillator.  The input function $I(t)$, the threshold
function $g(t)$ and the reset function $h(t)$ are deterministic
functions and will be regarded as given as part of the problem.
The parameter $\gamma \ge 0$ gives the rate of leakage of current
across the membrane.  The terms ``leaky'' and ``non-leaky'' are
sometimes used to describe the cases $\gamma > 0$ and $\gamma = 0$
respectively.  The issue is to provide a concise description of
the distribution of the random sequence $\{\tau_n^\e: n \ge 1\}$.
It is of particular interest to describe how the distribution of
the sequence $\{\tau_n^\e: n \ge 1\}$ responds to changes in one
or more of the functions $I$, $g$ and $h$.

If all three functions are constant, then the inter-spike
intervals $\tau_{n+1}^\e - \tau_n^\e$ form an independent,
identically distributed sequence of random variables.  A more
interesting situation occurs when one of the functions undergoes a
periodic modulation.  In many applications the input function is
taken to be of the form $I(t) = I_0 + I_1 \sin \omega t$. In other
cases the threshold is taken to be of the form $g(t) = g_0+g_1
\sin \omega t$. In this paper we will make the general assumption
that the three functions $I$ and $g$ and $h$ all have the same
period. Without loss of generality we will assume that the period
is 1. Then the sequence of firing phases
   $$
   \Theta_n^\e \equiv \tau_n^\e \quad \mbox{mod } 1
   $$
determines a Markov chain $\{\Theta_n^\e: n \ge 1\}$ on the circle
$\mathbb{S} = \mathbb{R}/\mathbb{Z}$.

When $\e = 0$ the process $X_t^0$ is given by an ordinary
differential equation.  Therefore the hitting times are given by
$\tau_{n+1}^0 = f(\tau_n^0)$ for some deterministic function $f$
satisfying $f(t+1) = f(t)+1$, and the firing phases are given by
$\Theta_{n+1}^0 = \tilde{f}(\Theta_n^0)$ where $\tilde{f}(\theta)
\equiv f(\theta)$ mod 1. In settings where either the input
function $I(t)$ or the threshold function $g(t)$ is of the form
$A+B \sin 2\pi t$ (and the other two functions are constant), the
dynamical system on $\mathbb{S}$ generated by iterating
$\tilde{f}$ has been studied by Rescigno, Stein, Purple and
Poppele \cite{rspp}, Knight \cite{kni}, Glass and Mackey \cite{gm}
and Keener, Hoppensteadt and Rinzel \cite{khr}. Of particular
interest are the regions in the $(A,B)$ parameter space giving
rise to phase-locked behavior, and the bifurcation scenario as $A$
and $B$ are varied.

When $\e > 0$ the deterministic hitting time function $f$ is
replaced by the first passage-time density function
$$
p^\e(t|t_0) := \frac{\partial}{\partial t}\mathbb{P}(\tau_n^\e
\leq t| \tau_{n-1}^\e = t_0)
$$
and $\tilde{f}$ is replaced the projection
$\tilde{p}(\theta|\theta_0)$, say, of $p^\e(t|t_0)$ onto the
circle $\mathbb{S}$. The behavior of the Markov chain
$\{\Theta_n^\e: n \ge 1\}$ may be studied via its transition
operator $T^\e$ given by
   $$
   T^\e \phi(\theta) = \mathbb{E}\left(\phi(\Theta_1^\e )\big| \Theta_0^\e =
   \theta\right) = \int_{\mathbb{S}}
\phi(\theta)\tilde{p}^\e(\theta|\theta_0)d\theta = \sum_m
\int_{\mathbb{S}} \phi(\theta)p^\e(\theta+m|\theta_0)d\theta
   $$ for $\phi$
in the class $B(\mathbb{S})$ of bounded measurable functions on
$\mathbb{S}$.   For any $\e > 0$ the transition densities
$\tilde{p}^\e(\theta|\theta_0)$ are bounded away from zero, so
that the Markov chain $\{\Theta^\e_n: n \ge 1\}$ is uniformly
ergodic and has a unique stationary probability distribution.  The
compact operator $T^\e$ captures the essential dynamics of
$\Theta_n^\e$, and hence its spectrum is of primary interest in
quantifying the transient and asymptotic behavior of the system.

In a sequence of papers Tateno \cite{t1,t2} and Tateno and Jimbo
\cite{tj} consider the effect of small noise on the deterministic
bifurcation scenarios considered earlier.  The papers
\cite{t1,t2,tj} contain numerical calculations of the leading
eigenvalues of the transition operator $T^\e$.  These calculations
suggest a qualitative change in the small noise behavior of the
leading eigenvalues near the location of the deterministic
bifurcation.  The calculations in \cite{t1,t2,tj} involve
numerical approximations in two places. Firstly, since there no
explicit formula for the first-passage density $p(t|t_0)$ except
in a few special cases, numerical techniques are used to solve an
integral equation for $p(t|t_0)$, following the method proposed by
Buonocore, Nobile and Ricciardi \cite{bnr}.  Secondly, the circle
$\mathbb{S}$ is replaced by a finite set of points. Thus the
operator $T^\e$ acting on $B(\mathbb{S})$ is approximated by a
finite-dimensional stochastic matrix.

In this paper we obtain rigorous results on the asymptotic
behavior of the spectrum of the operator $T^\e$ as $\e \to 0$. The
first main result gives a Gaussian approximation for the
first-passage density $p^\e(t|t_0)$ as $\e \to 0$.  This result
does not use the assumption of periodicity, and is valid for any
$C^2$ functions $I(t)$ and $g(t)$ and starting position $X^\e(t_0)
= x_0 < g(t_0)$, under the condition that the deterministic
hitting time is finite and that the deterministic trajectory
crosses the threshold transversally.  For details see Section
\ref{sec-passage} and especially Theorem \ref{thm-unif}.   This
result can be applied in the periodic setting to show that the
Markov chain $\{\Theta^\e_n: n \ge 1\}$ can be well approximated
by small Gaussian perturbations away from the deterministic
mapping $\tilde{f}$. The estimate in Theorem \ref{thm-unif} is
sufficiently strong that the techniques in Mayberry \cite{jm},
which deals with small Gaussian perturbations of circle maps, can
be applied here also.   In the simplest case where the
deterministic mapping $\tilde{f}$ is continuous and has one stable
fixed point $\theta_s$ attracting all orbits except the one
started at one unstable fixed point $\theta_u$, the limiting
eigenvalues of $T^\e$ can be calculated explicitly in terms of
$\tilde{f}'(\theta_s)$ and $\tilde{f}'(\theta_u)$, see Theorem
\ref{thm-outrans} in Section \ref{sec-sifmtrans}.  This result can
be extended to the case where $\tilde{f}$ is continuous and
phase-locked, see Remark \ref{rmk-periodp}.  In many examples of
SIF the deterministic mapping $\tilde{f}$ has a finite set of
discontinuities, and this case is treated in Section
\ref{sec-disc}.  The main result, Theorem \ref{thm-discspec},
deals with the case where $\tilde{f}$ is phase-locked and where
the discontinuities are well away from the phase locked orbit.
(This rather vague assertion is made precise in condition {\bf
(D3)} of Theorem \ref{thm-discspec}.)

Tateno and Jimbo \cite{tj} consider the leaky SIF with constant
input $I$, periodically modulated threshold $g(t) = 1+k\sin 2 \pi$
and constant reset level $0$. Using a $100 \times 100$ stochastic
matrix in place of the operator $T^\e$, they produce plots of the
leading eigenvalues for various small values of the noise
intensity $\e$.  In Section \ref{sec-examples}, we indicate how,
in the phase-locked setting, our results may be applied to give a
theoretical interpretation of the $\e \to 0$ behavior seen in some
of the figures of \cite{tj}.

Finally, Sections \ref{sec-passageproofs},
\ref{sec-sifmtransproofs} and \ref{sec-discproofs} contains the
proofs for the results in Sections \ref{sec-passage},
\ref{sec-sifmtrans} and \ref{sec-disc} respectively.

\section{First passage times} \label{sec-passage}

The results in this section do not use periodicity, and are valid
for general $C^2$ functions $I(t)$ and $g(t)$, and for all $\gamma
\ge 0$.

Denote by $\PP^{t_0,x_0}$ the law of the diffusion process
$\{X_t^\e: t \ge t_0\}$ satisfying
\begin{equation} \label{diff}
dX_t^\e =  (-\gamma X_t^\e +I(t))dt + \e \, dW_t
  \end{equation}
with initial condition $X_{t_0}^\e = x_0$.  For $x_0 < g(t_0)$
define the first passage time
   $$
   \tau^\e = \inf\{t \ge t_0: X_t^\e = g(t)\}
   $$
of the process $X_t^\e$ across the threshold $g(t)$.   In this
section we will consider the behavior of the distribution of
$\tau^\e$, and in particular its density function
  $$
  p^\e(t|t_0,x_0) = \frac{\partial}{\partial t}
  \PP^{t_0,x_0}(\tau^\e \le t),
  $$
as $\e \to 0$.

Let $\xi(t)=\xi(t|\, t_0,x_0)$ denote for $t \ge t_0$ the solution
to the noise free ($\e = 0$) equation (\ref{diff}) with initial
condition $\xi(t_0) = x_0 < g(t_0)$.   Thus
   $$
\xi(t|t_0,x_0) = e^{-\gamma(t-t_0)}x_0 +\int_{t_0}^t
e^{-\gamma(t-s)}I(s)\,ds.
   $$
For $x_0 < g(t_0)$ let
   $$
   f(t_0,x_0) = \inf\{t \ge t_0: \xi(t|t_0,x_0) = g(t)\}
   $$
denote the deterministic hitting time.  If $f(t_0,x_0) < \infty$,
define
  $$
   m(t_0,x_0) = -\gamma g(f(t_0,x_0))+ I(f(t_0,x_0))-g'(f(x_0,t_0)).
  $$
Thus $m(t_0,x_0)$ measures the difference in slopes when the
deterministic solution $\xi(t|t_0,x_0)$ first meets the threshold
$g(t)$.  Since $x_0 < g(t_0)$, the deterministic solution hits
from below and so $m(t_0,x_0) \ge 0$. The deterministic solution
crosses the threshold transversally if and only if $m(t_0,x_0) >
0$. Now define the set
$$
  {\cal G} = \{(t_0,x_0): x_0 < g(t_0) \mbox{ and }
f(t_0,x_0) < \infty \mbox{ and } m(t_0,x_0) > 0\}.
  $$
of initial conditions $(t_0,x_0)$ for which the deterministic
trajectory crosses the threshold function transversally at the
finite time $f(t_0,x_0)$.

\begin{proposition} \label{f}
${\cal G}$ is an open set and $f \in C^2({\cal G},\mathbb{R})$.
\end{proposition}

The proofs of the results in this section can be found in Section
\ref{sec-passageproofs}.  Our next result gives a uniform bound on
the deviation of $\tau^\e$ away from the deterministic crossing
time $f(t_0,x_0)$.

\begin{proposition} \label{deviations}
Let $G$ be a compact subset of ${\cal G}$. Then for any $\delta
>0$ there are constants $M_\delta$ and $K_\delta$ such that
\begin{equation*}
\mathbb{P}^{t_0,x_0}(|\tau^\e-f(t_0,x_0)| > \delta) \leq M_\delta
\e \, e^{-K_\delta/\e^2}
\end{equation*}
for all $(t_0,x_0) \in G$.
\end{proposition}

\s

For $(t_0,x_0) \in  {\cal G}$ define
     $$
\sigma^2(t_0,x_0) =
\left\{\begin{array}{cl}(1-e^{-2\gamma(f(t_0,x_0)-t_0)})/2\gamma &
\hspace{.3in} \mbox{ if }\gamma > 0, \\[1ex]
 f(t_0,x_0)-t_0 &
\hspace{.3in} \mbox{ if }\gamma = 0. \end{array}\right.
   $$
Notice that $\e^2 \sigma^2(t_0,x_0)$ is the variance of $X_t^\e$
at the moment $t = f(t_0,x_0)$ of noise-free intersection.  Define
    \begin{equation} \label{sigma_tau}
\sigma_\tau^2(t_0,x_0) =  \frac{\sigma^2(t_0,x_0)}{m^2(t_0,x_0)}
    \end{equation}
and
 $$
 p_\tau(t|t_0,x_0) = \frac{1}{\sqrt{2\pi\sigma^2_\tau(t_0,x_0)}}
 e^{-t^2/2\sigma^2_\tau(t_0,x_0)},
 $$
so that $p_\tau(t|t_0,x_0)$ is the density at $t$ of a
$N(0,\sigma^2_\tau(t_0,x_0))$ normal random variable. With this
notation in hand, we can state the main result of this section.

\begin{thm} \label{thm-unif} Let $G$ be a compact subset of ${\cal G}$. Then there exist finite positive constants $\delta$, $\sigma_1$,
$K$ and $\e_0$ (depending on $G$) so that
\begin{equation} \label{clt}
 \sup_{(t_0,x_0) \in G} |\e p^\e(f(t_0,x_0)+u|t_0, x_0)- p_\tau(u/\e|t_0,x_0)| \leq K \e e^{-u^2/2\e^2\sigma_1^2}
\end{equation}
for all $\e < \e_0$, $|u| \leq \delta$.
\end{thm}

\begin{figure}[ht]
\begin{center}
\includegraphics[scale=0.7]{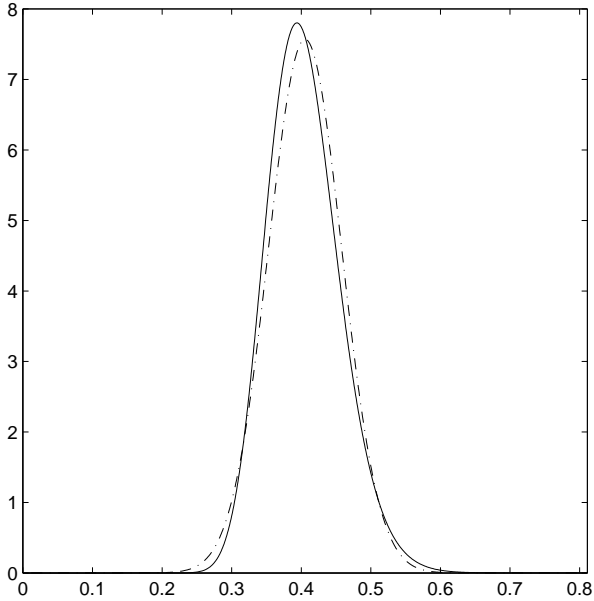}
\end{center}
\caption {First passage-time pdf for $dX_t = (-X_t+2)dt+.1dW_t$,
$X_0=.5$, across $g(t)=1$. Dashed line: Gaussian approximation of
Theorem \ref{thm-unif} with mean $\ln(1.5) \approx .4055$ and
standard deviation $\sqrt{5/1800} \approx 0.0527$. Solid:
Numerical approximation obtained by solving integral equation from
\cite{bnr}} \label{passagefig}
\end{figure}

\begin{remark} \label{rmk-gaussian}
{\rm Under $\PP^{t_0,x_0}$ the centered and scaled hitting time
$(\tau^\e-f(t_0,x_0))/\e$ has density $\e p^\e(f(t_0,x_0)+\e
t|t_0,x_0)$ at $t$.    Putting $u  =\e t$ in Theorem
\ref{thm-unif} gives the result that under $\PP^{t_0,x_0}$
   \begin{equation} \label{cltdist}
    \frac{\tau^\e-f(t_0,x_0)}{\e} \Rightarrow
    N(0,\sigma^2_\tau(t_0,x_0)) \quad \mbox{ as } \e \to 0,
   \end{equation}
or more informally
    $$
    \tau^\e \approx N(f(t_0,x_0), \e^2 \sigma^2_\tau(t_0,x_0)) \quad \mbox{ as } \e \to
    0.
    $$
Figure \ref{passagefig} shows the densities of $\tau^\e$ and
$N(f(t_0,x_0), \e^2 \sigma_\tau(t_0,x_0))$ for an example with $\e
= 0.1$.  Of course the result in Theorem \ref{thm-unif} is much
stronger than (\ref{cltdist}) since it gives locally uniform
convergence of densities, rather than just convergence in
distribution.  This extra strength will be important for the
spectral analysis results in the next two sections.}
\end{remark}

\begin{remark} \label{rmk-slope}
{\rm In equation (\ref{sigma_tau}) the difference in slopes
$m(t_0,x_0)$ at the point of deterministic intersection is used to
convert the variance $\sigma(t_0,x_0)$ in the spatial dimension
into the variance $\sigma_\tau(t_0,x_0)$ in the temporal
dimension.  A heuristic observation of this conversion factor
appears in Stein \cite{ste}.}
\end{remark}

\section{Transition Operator for the SIF} \label{sec-sifmtrans}

   We now return to the setting of SIFs.
The functions $I(t)$, $g(t)$ and $h(t)$ are assumed to be $C^2$
and periodic with period 1. Also $h(t) < g(t)$ for all $t$, and
$\gamma \ge 0$.  The deterministic solution is
   $$
   \xi(t|t_0,x_0) = e^{-\gamma(t-t_0)}x_0 + \int_{t_0}^t e^{-\gamma(t-s)}
   I(s)\,ds.
   $$
If $\gamma = 0$, then the assumption $\int_0^1 I(s)\,ds > 0$
implies that $\xi(t|t_0,x_0) \to \infty$ as $t \to \infty$ and so
the deterministic hitting $f(t_0,x_0)$ is finite.  If $\gamma > 0$
then $\left|\xi(t|t_0,x_0) - \overline{\xi}(t)\right| \to 0$ as $t
\to \infty$ where
   $$
   \overline{\xi}(t) = \int_{-\infty}^t e^{-\gamma(t-s)}
   I(s)\,ds = \frac{1}{1-e^{-\gamma}}\int_{t-1}^t e^{-\gamma(t-s)}
   I(s)\,ds = \frac{1}{e^\gamma - 1} \int_0^1 e^{\gamma u}I(t+u)\,du.
   $$
In order to ensure that crossings for the noise free system occur
in a nice enough fashion we impose the following conditions on the
functions $I(t)$ and $g(t)$:
  \begin{description}
  \item[(A)] $ \left\{ \begin{array}{cl}
    \max\{\overline{\xi}(t) - g(t): 0 \le t \le 1\} > 0 & \mbox{ if } \gamma > 0
    \\[1ex]
    \displaystyle{\int_0^1 I(s)ds > 0 }& \mbox{ if } \gamma = 0
      \end{array} \right.  $
  \item[(B)] $-\gamma g(t) + I(t) - g'(t) > 0$ for all $t$
  \end{description}
Note that {\bf (A)} implies that $f(t_0,x_0) < \infty$ whenever
$x_0 < g(t_0)$, and then {\bf (B)} implies that $m(t_0,x_0) > 0$.
Therefore ${\cal G} = \{(t_0,x_0): x_0 < g(t_0)\}$, and we can
apply our results on first passage densities from Section
\ref{sec-passage} using the compact set $G = \{(s,h(s)): 0 \le s
\le 1\}$.  We always take $x_0 = h(t_0)$, so we write
$p^\e(t|t_0)=p^\e(t|t_0,h(t_0)) $, $f(t_0) = f(t_0,h(t_0))$,
$\sigma_\tau^2(t_0)
 = \sigma^2_\tau(t_0,h(t_0))$, etc.\
Proposition \ref{f} implies that $f \in C^2(\R)$, and clearly
$f(t+1)=f(t) + 1$ for all $t \in \mathbb{R}$.  Theorem
\ref{thm-unif} implies that when $\e > 0$ is small, we have the
approximation
 \begin{equation} \label{pertapp}
 \tau_n^\e \approx f(\tau_{n-1}^\e) + \e \sigma_\tau(\tau_{n-1}^\e)
 \chi_n
 \end{equation}
where $\{\chi_n: n \ge 1\}$ is a sequence of independent $N(0,1)$
random variables.

The sequence of firing phases $\{\Theta_n^\e: n \ge 1\}$ is a
Markov chain on the circle $\mathbb{S} = \mathbb{R}/\mathbb{Z}$
with transition density function
\begin{equation} \label{transdensity}
\tilde{p}^\e(\theta|\theta_0) = \sum_{m \in \, \mathbb{Z}}
p^\e(\theta + m|\theta_0)
\end{equation}
for all $\theta,\theta_0 \in \mathbb{S}$.  Looking at
\eqref{pertapp}, we may expect that the study of $\{\Theta_n^\e: n
\ge 1\}$ should be similar to the study of the chain
\begin{equation} \label{MC}
Y_n^\e = f(Y_{n-1}^\e) + \e \sigma_\tau(Y_{n-1}^\e) \chi_n
\hspace{.2in} \mod 1.
\end{equation}
Spectral properties of the transition operator for Markov chains
of the form \eqref{MC} were developed in Mayberry \cite{jm}. Here,
we prove a similar result for the transition operator $T^\e$ of
the chain $\{\Theta_n^\e: n \ge 1\}$.

\begin{thm} \label{thm-outrans}
Let $\{\Theta_n^\e: n \ge 1\}$ be the sequence of firing phases
for a period 1 SIF with $C^2$ input, threshold and reset functions
$I(t)$, $g(t)$ and $h(t)$.  Assume the conditions {\bf(A)} and
{\bf(B)}. Suppose that the deterministic phase return map
$$
\tilde{f} \equiv f \quad \mod \, 1
$$
has a stable fixed point $\theta_s$ and an unstable fixed point
$\theta_u$, and that $\tilde{f}^n(\theta) \rightarrow \theta_s$
for all $\theta \in \mathbb{S} \setminus \{\theta_u\}$.  Let
$T^\e$ denote the transition operator for $\{\Theta_n^\e: n \ge
1\}$. Then for any $r
> 0$ and $\e$ sufficiently small, we can write $T^\e = T_{lp}^\e +
T_{up}^\e$ where $\|T_{lp}^\e\|_\infty < r$ and any eigenvalue of
$T_{up}^\e$ with modulus greater than $r$ is of one of the two
forms $c_s^n + O(\e)$ or $|c_u|^{-1}c_u^{-n} + O(\e)$ for some $n
\geq 0$, where $c_s = f'(\theta_s)$ and $c_u = f'(\theta_u)$.
\end{thm}

The proof of this result can be found in Section
\ref{sec-sifmtransproofs}.  The result says that $c_s^n$ and
$|c_u|^{-1}c_u^{-n}$ are the {\bf limiting eigenvalues} of $T^\e$
in the sense that for any $r
> 0$, $T^\e$ has sequences of $r$-pseudoeigenvalues which converge
to $c_s^n$ and $|c_u|^{-1}c_u^{-n}$  as $\e \to 0$ (see Trefethen
and Embree \cite{te} for definitions of pseudoeigenvalues).

\begin{remark} \label{rmk-periodp}
{\rm Theorem \ref{thm-outrans} can be extended to the case where
$\tilde{f}$ has periodic orbits.   Consider an orbit
$\{\theta_1,\theta_2, \ldots,\theta_\kappa\}$ of period $\kappa \ge 1$ and
let $c = f'(\theta_1)f'(\theta_2) \cdots f'(\theta_\kappa)$ denote the
product of the derivatives of $f$ along the periodic orbit. If the
orbit is stable, so that $|c| < 1$, then it contributes limiting
eigenvalues $(c^n)^{1/\kappa}$ for $n \ge 0$.  If the orbit is
unstable, so that $|c| > 1$, then it contributes limiting
eigenvalues $(|c|^{-1}c^{-n})^{1/\kappa}$.  Here all the $\kappa$th roots
are included as limiting eigenvalues.  The proof in this more
general setting combines the method of proof of Theorem
\ref{thm-outrans} with techniques used in the proof of
\cite[Theorem 1]{jm}, and is left to the reader. Moreover, the
methods used in \cite[Theorem 2]{jm} can be applied here to give
information about the associated eigenfunctions.}
\end{remark}

\section{Discontinuous Case} \label{sec-disc}

Tateno and Jimbo \cite{tj} consider several cases of SIFs in which
$f$ is well defined, but Assumption {\bf (B)} fails and $f$ is
discontinuous at some $\theta^* \in \mathbb{S}$. In this section,
we will discuss extensions of Theorem \ref{thm-outrans} to this
situation. As before, for $h(t_0) < g(t_0)$ we define $f(t_0) =
\inf\{t \geq t_0 : \xi(t|t_0) = g(t) \} < \infty$, but now we also
define $f^*(t_0) = \inf\{t > t_0: \xi(t|t_0) > g(t)\}$.  Thus
$f(t_0)$ is the time of first hitting of the threshold, and
$f^*(t_0)$ is the time of first crossing of the threshold.  We
keep {\bf (A)} unchanged, but replace condition {\bf (B)} with
  \begin{description}
   \item[(B')] There is a finite set $D \subset [0,1)$ (possibly empty) such that $-\gamma g(f(t_0)) + I(f(t_0))- g'(f(t_0))> 0$ for all $t_0 \in [0,1)\setminus D$.
   \item[(C')]  For each $t_0 \in D$ either
     \begin{description}
      \item[(i)] $f^*(t_0) = f(t_0)$; or else
      \item[(ii)] $f^*(t_0) > f(t_0)$ and $g(t) > \xi(t|t_0)$ for $f(t_0) < t < f^*(t_0)$.
   \end{description}
   \end{description}

\begin{figure}[h]
\begin{center}
\includegraphics[scale=0.75]{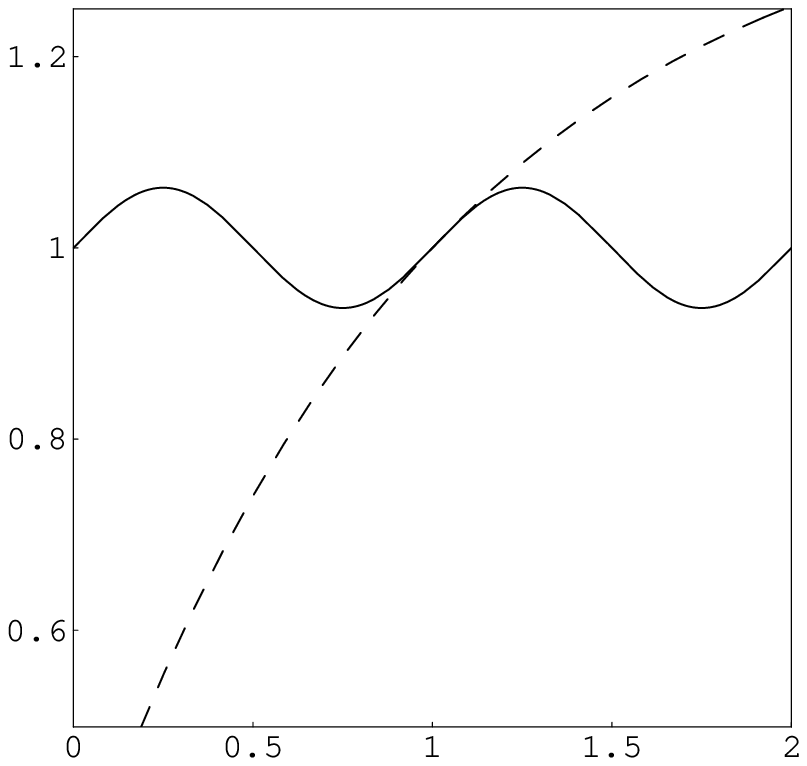} \hspace{6ex}
\includegraphics[scale=0.75]{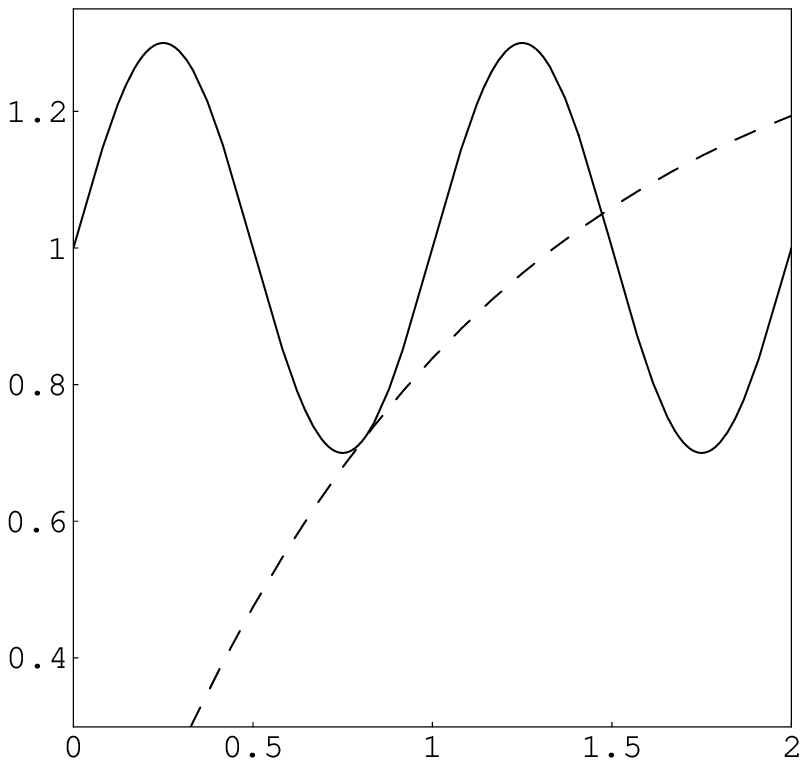}
\end{center}
\caption{Examples of the behavior described in {\bf (C')}.  In
both cases $\gamma=1$, $I(t)\equiv 1.4$, $g(t)=1+B \sin 2\pi t$
and $h(t)\equiv 0$.  Case {\bf(i)} on the left has $B=.0629,
t_0=-.2527$ giving $f(t_0)=f^*(t_0)=1.0251$.  Case {\bf (ii)} on
the right has $B=.3, t_0=.0863$ giving $f(t_0)=.8087$,
$f^*(t_0)=1.473$.} \label{Ccond}
\end{figure}

Notice that {\bf (A)} implies that $f(t_0) \le f^*(t_0) < \infty$
for all $t_0$, and that {\bf (B')} implies that Proposition
\ref{deviations} can be applied to any compact subset of $[0,1]
\setminus D$.  In case {\bf (C')(i)} the deterministic trajectory
$\xi(t|t_0)$ crosses the threshold $g(t)$ at $t=f(t_0)$, and $f$
is continuous (but not differentiable) at $t_0$; and in case {\bf
(C')(ii)} $\xi(t|t_0)$ touches $g(t)$ at $t=f(t_0)$ and then does
not intersect again until it crosses at time $t = f^*(t_0)$, and
$f$ is discontinuous at $t_0$.  Examples of the behavior described
in {\bf (C')} are given in Figure \ref{Ccond}.

\s

We begin with the following extension of Proposition \ref{deviations}. Proofs for the results in this section can be found in Section \ref{sec-discproofs}.

\begin{proposition} \label{deviations2}
Suppose that $f$ satisfies {\bf (A)}, {\bf (B')} and {\bf (C')},
and $D \neq \emptyset$.  Then for all $\delta
>0$ there exist $\tilde{\delta} > 0$ and $K,M$ such that
  $$
  \PP^{t}(d(\tau^\e,\{f(t_0),f^*(t_0)\}) > \delta) \leq \e M
e^{-K/\e^2}
$$
whenever $|t-t_0| < \tilde{\delta}$ for some $t_0 \in D$.
\end{proposition}

\s

The assumptions in our final theorem are not the most general ones
possible, but the result is sufficient to treat the examples in
the next section and to give the reader the indication of how
Theorem \ref{thm-outrans} (see also Remark \ref{rmk-periodp}) can
be extended to discontinuous settings. The conditions may seem
awkward, but we will see in the next section that they are easy to
verify numerically in examples
 of interest.

\begin{thm} \label{thm-discspec}
Suppose that {\bf(A)}, {\bf(B')}, and {\bf(C')} are satisfied and
in addition that $\tilde{f}=f$ mod 1 satisfies the conditions
\begin{description}
\item[(D1)] $\tilde{f}$ has a periodic orbit $P = \{\theta_1, \ldots,\theta_\kappa\}$ in $\mathbb{S} \setminus D$ for some $\kappa \ge 1$
and $$|\tilde{f}'(\theta_1)\tilde{f}'(\theta_2) \cdots
\tilde{f}'(\theta_\kappa)| < 1.$$
\item[(D2)] $\tilde{f}^n(\theta) \to P$ as $n \to \infty$ for all $\theta \in \mathbb{S}$.
\item[(D3)] $\tilde{f}^{-\ell}(D) = \emptyset$ for some $\ell \ge 1$ and $\tilde{f}^i(E) \cap D = \emptyset$
  for $0 \le i \le \ell - 1$, where $E = \tilde{f}(D) \cup \tilde{f}^*(D)$.
\end{description}
Let $T^\e$ denote the transition operator for $\Theta_n^\e$. Then
for any $r > 0$ and $\e$ sufficiently small, we can write $T^\e =
T_{lp}^\e + T_{up}^\e$ where $\|T_{lp}^\e\|_\infty < r$ and any
eigenvalue of $T_{up}^\e$ with modulus greater than $r$ is of the
form $(c^n)^{1/\kappa} + O(\e)$ for $n \ge 0$, where $c =
\tilde{f}'(\theta_1)\tilde{f}'(\theta_2) \cdots
\tilde{f}'(\theta_\kappa) = (\tilde{f}^\kappa)'(\theta_i)$ for any $i \in
\{1,\ldots,\kappa\}$.
\end{thm}

\section{Examples from Tateno and Jimbo} \label{sec-examples}

We now apply our spectral results to the examples considered in
\cite{tj}. There, the authors consider leaky SIFs with constant
input $I(t) \equiv I$, threshold $g(t) = 1 + k \sin 2\pi t$, and
reset $h(t)\equiv 0$.  Since $\overline{\xi}(t) = I/\gamma$ the
condition {\bf (A)} becomes
\begin{equation} \label{finite}
  I/\gamma > 1-k.
\end{equation}
The transversality condition {\bf (B)} is now
$$
-\gamma( 1 +k \sin 2\pi t) + I - 2 \pi k \cos 2 \pi t > 0  \quad
\quad \mbox{ for all }t
$$
or equivalently
\begin{equation} \label{transversal}
I/\gamma - 1 > k \sqrt{4\pi^2/\gamma^2 + 1}.
\end{equation}
If condition {\bf (B)} fails, then there is a point of tangency on
the threshold curve which leads to a point with the property
described in {\bf (C')(i)} or {\bf (C')(ii)}.  In any case,
assuming {\bf (A)} holds, we obtain $f(t_0) = \inf\{t \geq t_0:
h(t) = t_0\}$ where
  $$
  h(t) = t+ \frac{1}{\gamma} \log\left(1 - \frac{\gamma}{I}(1+k\sin 2\pi
  t)\right).
  $$
If $h'(t) > 0$ for all $t$ then $f = h^{-1}$ is a smooth function.
However, strict local maxima of $h$ give rise to discontinuities in
$f$.  Moreover periodic orbits of the induced mapping $\tilde{f}$ on
$\mathbb{S}$ can be found using the facts that if $h^\kappa(t_0) =
t_0$ mod 1 then $f^\kappa(t_0) = t_0$ mod 1 and $(f^\kappa)'(t_0) =
1/(h^\kappa)'(t_0)$.  The plots of $\tilde{f}$ and the numbers in
the examples below are obtained using the explicit formula for $h$.

\s

Tateno and Jimbo take $\gamma = 1/12.8 $ and various values of $I$
and $K$.

\s

\begin{figure}[h]
\begin{center}
\includegraphics[scale=0.8]{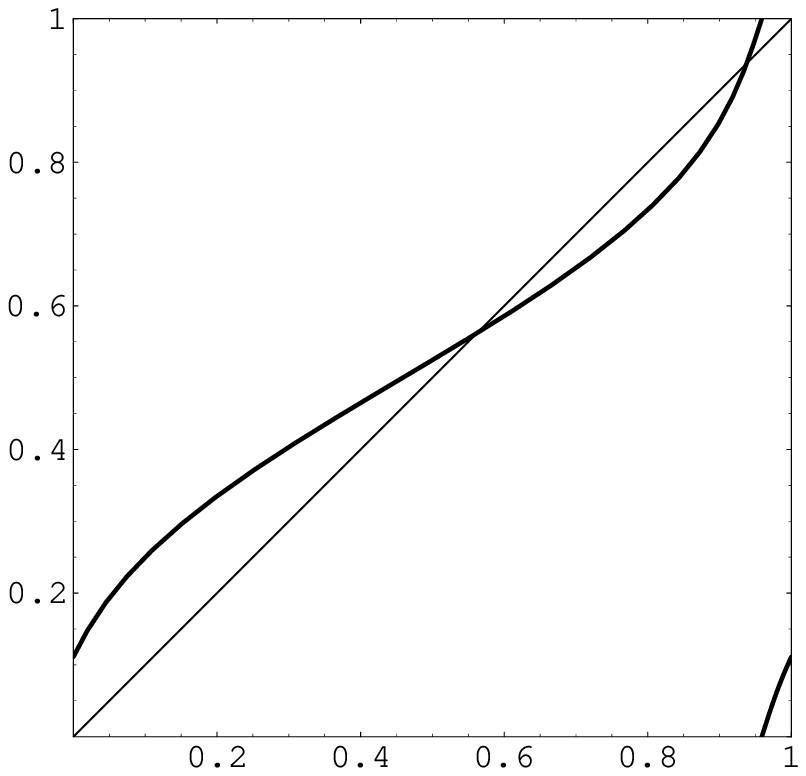}
\end{center}
\caption{Plot of $\tilde{f}$ for $I=1, k=0.1, \gamma=1/12.8$}
\label{I1kpt1}
\end{figure}

\n{\bf Example 1.} Taking $I = 1$ and $k = 0.1$ gives the map
$\tilde{f}$ in Figure \ref{I1kpt1}, which satisfies the conditions
of Theorem \ref{thm-outrans}. There is a stable fixed point at
$0.5622$ with $\tilde{f}'(0.5622) = 0.6142$ and an unstable fixed
point at $0.9379$ with $\tilde{f}'(0.9379) = 2.6898$.  All orbits
starting at $\theta_0 \neq 0.9379$ converge to $0.5622$ so by
Theorem \ref{thm-outrans} the limiting eigenvalues are $\{0.6142^n:
n \ge 0\} \cup \{2.6898^{-n-1}: n \ge 0\} = \{ 1, 0.6142, 0.3772,
0.3718, 0.2317, \ldots\}$

\s

\begin{figure}[ht]
\begin{center}
\includegraphics[scale=0.8]{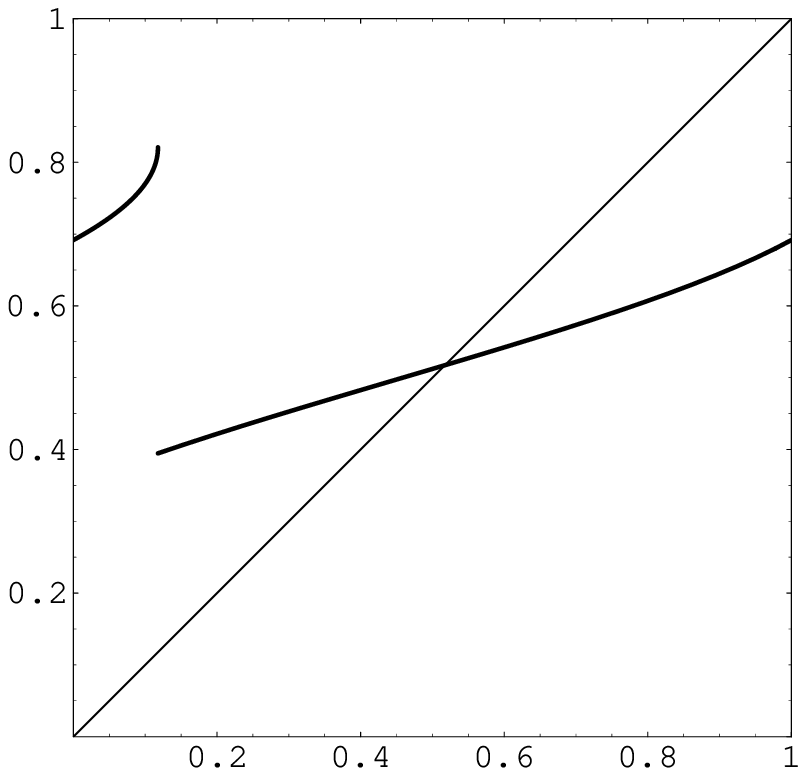}
\end{center}
\caption{Plot of $\tilde{f}$ for $I=1, k=0.35, \gamma=1/12.8$}
\label{I1kpt35}
\end{figure}

\n{\bf Example 2.}  Keeping $I = 1$ and increasing $k$ to $k = 0.35$
gives the map $\tilde{f}$ in Figure \ref{I1kpt35}.  There is a
discontinuity at 0.1178 with $\tilde{f}(0.1178) = 0.8208$ and
$\tilde{f}^*(0.1178) = 0.3946$ so in the notation of Theorem
\ref{thm-discspec}, we have $D = \{0.1178\}$ and $E =
\{0.8208,0.3946\}$.  Clearly $\tilde{f}^{-1}(D) = \emptyset$ and $E
\cap D = \emptyset$ so that {\bf(D3)} is satisfied with $\ell=1$.
There is a stable fixed point at 0.5173 with $f'(0.5173) = 0.2973$
so that {\bf(D1)} is satisfied with $\kappa=1$, and clearly {\bf
(D2)} is also satisfied. Therefore, by Theorem \ref{thm-discspec}
the limiting eigenvalues are $\{0.2973^n: n \ge 0\} = \{1, 0.2973,
0.0884, 0.0263, \ldots\}$.

\s

\s

The next four examples correspond closely to the values considered
by Tateno and Jimbo.

\s
\begin{figure}[htb]
\begin{center}
\includegraphics[scale=.8]{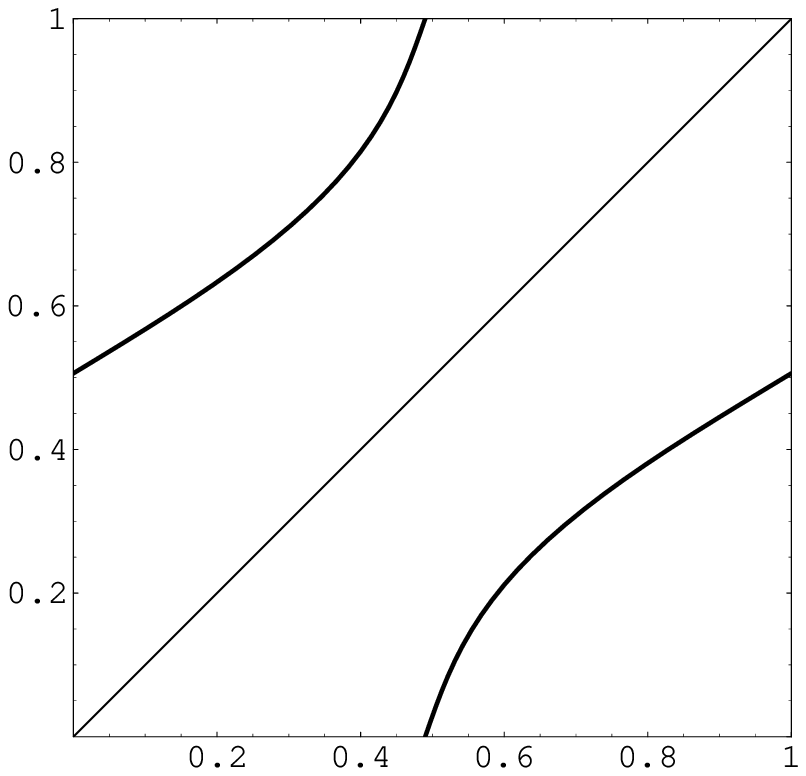} \hspace{6ex}
\includegraphics[scale=.8]{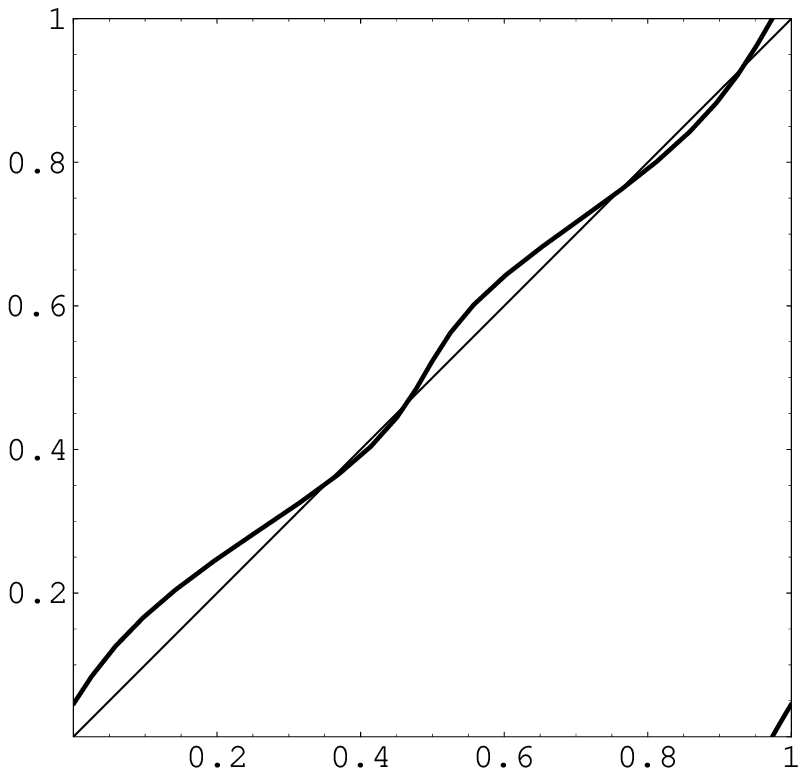}
\end{center}
\caption{Plot of $\tilde{f}$ (left) and $\tilde{f}^2$ (right) for
$I=2, k=0.2, \gamma=1/12.8$} \label{I2kpt2}
\end{figure}

\n{\bf Example 3.} $I = 2$ and $k = 0.2$.  Figure \ref{I2kpt2} shows
$\tilde{f}$ and $\tilde{f}^2$. There is a period 2 stable orbit
$\{0.3527, 0.7593\}$ with $(f^2)'(0.3527) = 0.7445$ and a period 2
unstable orbit $\{0.4654, 0.9329\}$ with $(f^2)'(0.4654) =1.5043$.
By Remark \ref{rmk-periodp} following Theorem \ref{thm-outrans} the
limiting eigenvalues are $\{ (0.7445^n)^{1/2}: n \ge 0\} \cup\{
(1.5043^{-n-1})^{1/2}: n \ge 0\} = \{\pm 1, \pm 0.8628, \pm 0.8153,
\pm 0.7445,  \ldots\}$.

\s

\begin{figure}[hbt]
\begin{center}
\includegraphics[scale=.8]{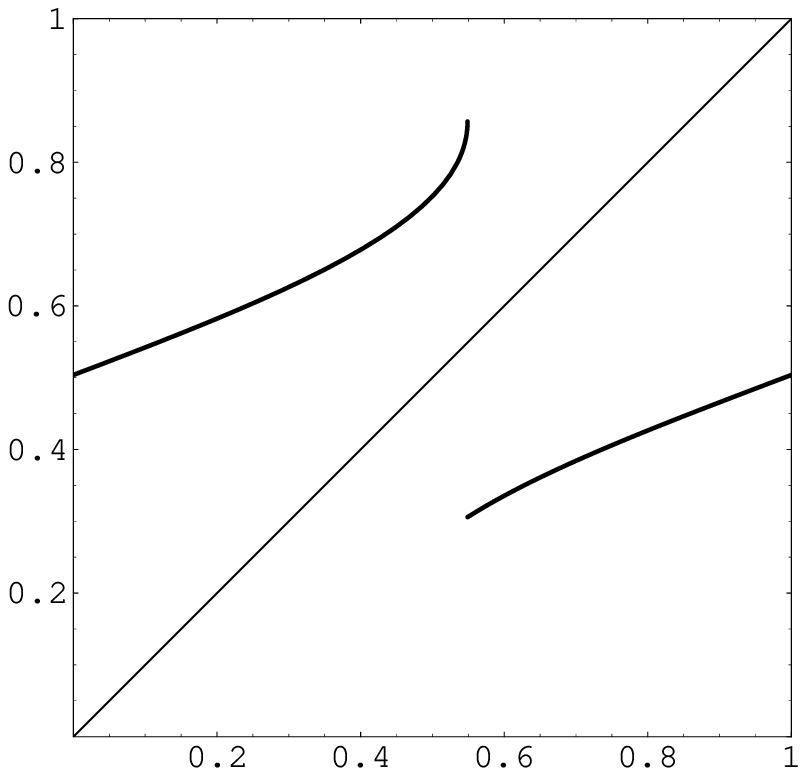} \hspace{6ex}
\includegraphics[scale=.8]{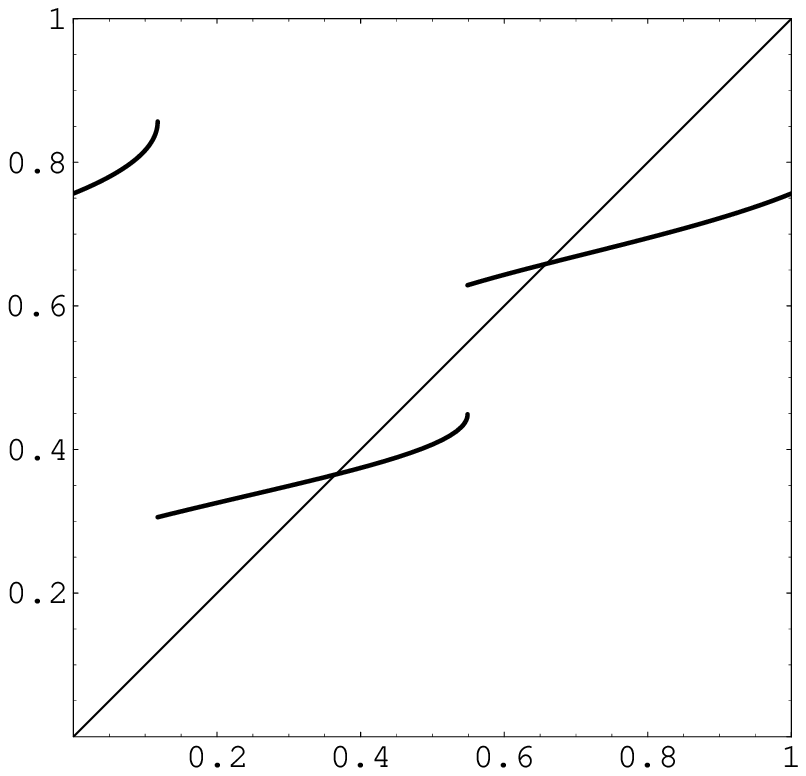}
\end{center}
\caption{Plot of $\tilde{f}$ (left) and $\tilde{f}^2$ (right) for
$I=2, k=0.5, \gamma=1/12.8$} \label{I2kpt5}
\end{figure}

\n{\bf Example 4.} $I = 2$ and $k = 0.5$.  Again we show $\tilde{f}$
and $\tilde{f}^2$, see Figure \ref{I2kpt5}.  There is a point of
discontinuity with $D = \{0.5489\}$ with $E= \{0.8567, 0.3057\}$. In
this example $\tilde{f}^{-1}(D) = \{.1174\} \neq \emptyset$, but
$\tilde{f}^{-2}(D) = \tilde{f}^{-1}(\{0.1174\}) = \emptyset$, and
{\bf (D3)} holds with $\ell=2$.  There is a period 2 stable orbit
$\{0.3651, 0.6586\}$, and the product of the derivatives along the
orbit is $(f^2)'(0.3651)= 0.2544$ %^{1/2} = 0.5043$
so that {\bf (D1)} holds with $\kappa=2$. It is clear from the plot
of $\tilde{f}^2$ that {\bf (D2)} also holds. Thus Theorem
\ref{thm-outrans} implies the limiting eigenvalues are $\{\pm
(0.2554^n)^{1/2}: n \ge 0\} = \{\pm 1, \pm 0.5044, \pm 0.2544,
\ldots\}$. This is one case considered in Figure 4 of Tateno and
Jimbo \cite{tj} and our predicted limiting values can be seen at the
extreme left edge of the $k = 0.5$ part of Figure 4 of \cite{tj}.

\s
\begin{figure}[ht]
\begin{center}
\includegraphics[scale=0.8]{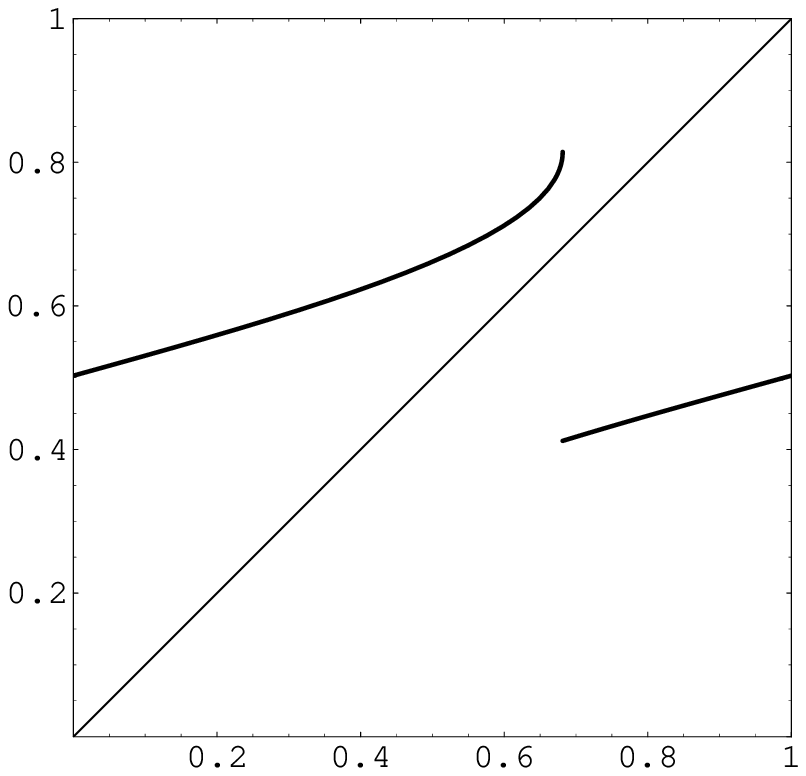} \hspace{6ex} \includegraphics[scale=0.8]{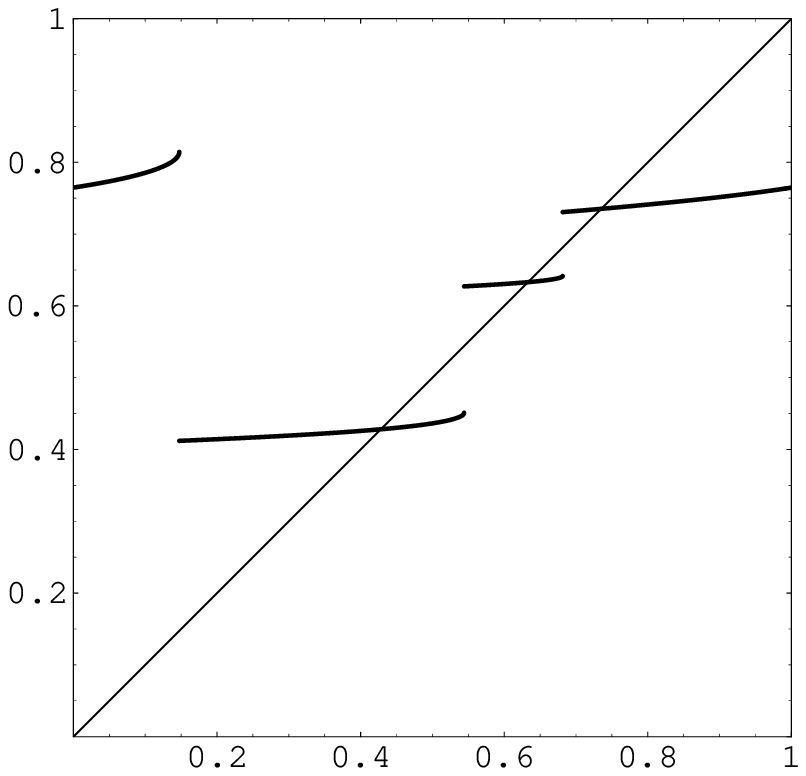}
\end{center}
\caption{Plot of $\tilde{f}$ (left) and third iterate
$\tilde{f}^3$ (right) for $I=2, k=0.8, \gamma=1/12.8$}
\label{I2kpt8}
\end{figure}

\n{\bf Example 5.} $I = 2$ and $k = 0.8$.  Figure \ref{I2kpt8}
shows $\tilde{f}$ and $\tilde{f}^3$.  We are again in the setting
of Theorem \ref{thm-discspec}. There is an attracting period 3
orbit $$\{0.4218,0.6330, 0.7352\}$$
and the product of the
derivatives along the orbit is $(\tilde{f}^3)'(0.4218) = 0.088076$
so the limiting eigenvalues are $\{(0.088076^n)^{1/3} : n \ge 0\}
= \{\omega^r 0.4449^n: r=0,1,2 \mbox{ and } n \ge 0\}$ where
$\omega = e^{2 \pi i/3}$ is a cube root of unity.  These can be
seen at the extreme left edge of the $k = 0.8$ part of Figure 4 of
\cite{tj}, and also in Figure 5(c) of \cite{tj}.

\s

\begin{figure}[h]
\begin{center}
\includegraphics[scale=0.8]{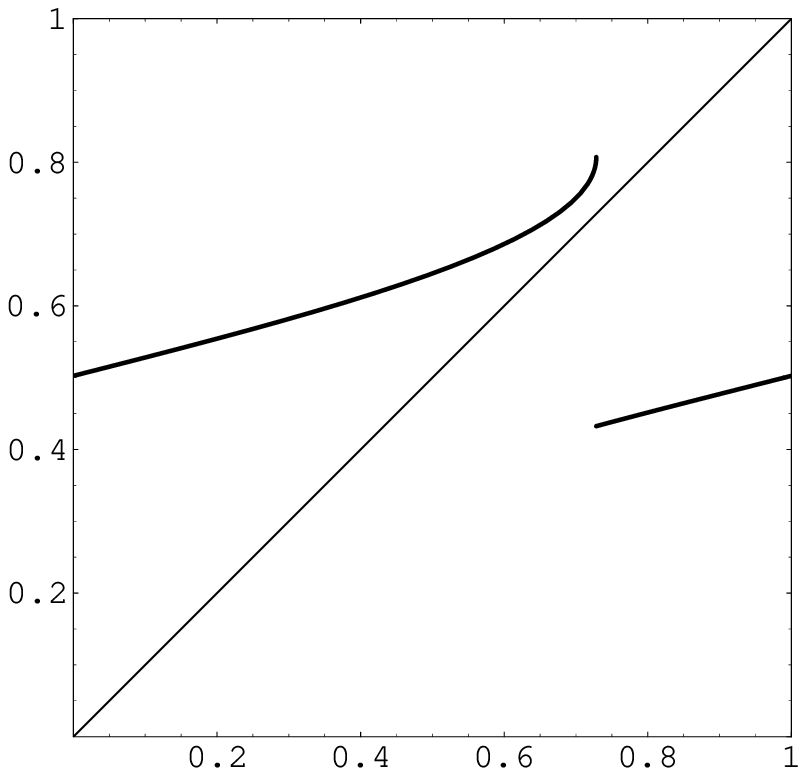} \hspace{6ex} \includegraphics[scale=0.8]{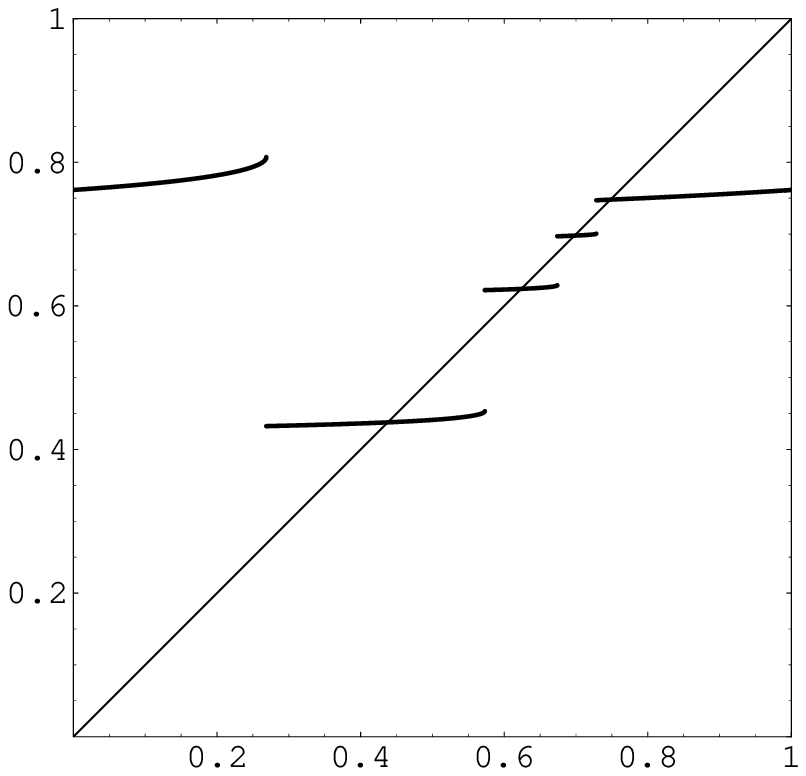}
\end{center}
\caption{Plot of $\tilde{f}$ (left) and fourth iterate
$\tilde{f}^4$ (right) for $I=2, k=0.9, \gamma=1/12.8$}
\label{I2kpt9}
\end{figure}
\n{\bf Example 6.} $I = 2$ and $k = 0.9$.  Figure \ref{I2kpt9}
shows $\tilde{f}$ and $\tilde{f}^4$. There is now an attracting
period four orbit $\{0.4378,0.6236, 0.6978, 0.7480\}$ and the
product of the derivatives along the orbit is $0.043991$ and the
conditions of Theorem \ref{thm-discspec} are satisfied so the
limiting eigenvalues are $\{(0.043991^n)^{1/4}: n \ge 0\} = \{i^r
0.4580^n: r=0,1,2,3 \mbox{ and } n \ge 0\}$.  This gives a
theoretical justification for the $k = 0.9$ part of Figure 5(c) of
\cite{tj} where fourth roots of one appear in the limiting
spectrum.

\section{Proofs for Section \ref{sec-passage}} \label{sec-passageproofs}

We show first that it suffices to prove the results for the case
of constant input function.  For any constant $I$ define
   $$
   k(t) = e^{-\gamma t} \int_0^t e^{\gamma s} (I(s)-I)\,ds .
   $$
Then $\widehat{X}^\e_t \equiv X^\e_t - k(t)$ satisfies
     $$
     d\widehat{X}^\e_t = \left(- \gamma \widehat{X}^\e_t + I\right)dt+\e
     dW_t.
     $$
Define $\widehat{g}(t) = g(t)-k(t)$, then $\widehat{X}^\e_t$
started at $x_0-k(t_0)$ at time $t_0$ hits the threshold
$\widehat{g}(t)$ at the same moment that $X^\e_t$ started at $x_0$
at time $t_0$ hits the threshold $g(t)$.   Moreover if
$\widehat\xi$, $\widehat{f}$, $\widehat{m}$, $\widehat{\cal G}$,
$\widehat{\sigma}$, $\widehat{\tau}^e$ and $\widehat{\PP}$ are
defined using the constant input $\widehat{I}(t)\equiv I$ and the
threshold function $\widehat{g}$ in the same way as $\xi$, $f$,
$m$, ${\cal G}$, $\sigma$, $\tau^e$ and $\PP $ are defined using
the original input $I(t)$ and original threshold $g(t)$, then
$\widehat{\xi}(t|t_0,x_0-k(t_0)) = \xi(t|t_0,x_0)-k(t)$,
$\widehat{f}(t_0,x_0-k(t_0)) = f(t_0,x_0)$,
$\widehat{m}(t_0,x_0-k(t_0)) = m(t_0,x_0)$, $\widehat{\cal G} =
\{(t_0,x_0-k(t_0)): (t_0,x_0) \in {\cal G}\}$,
$\widehat{\sigma}(t_0,x_0-k(t_0)) = \sigma(t_0,x_0)$ and
$\PP^{t_0,x_0-k(t_0)}(\widehat{\tau}^\e \in A) =
\PP^{t_0,x_0}(\tau^\e \in A)$ for any Borel subset $A \subset \R$.
 Therefore any of the results of Section \ref{sec-passage} proved under the
assumption that $I(t) \equiv I$ can be converted into the
corresponding result for a more general function $I(t)$.

This method of converting a problem with a time varying input into
one with a constant input and time varying threshold (and reset)
is used in Scharstein \cite{sch}.  For the reminder of this
section we shall assume that $I(t) \equiv I$, so that $X^\e_t$ is
either an Ornstein-Uhlenbeck process (if $\gamma
> 0$ or else a Brownian motion with constant drift.

\s

\n {\bf Proof of Proposition \ref{f}.}  Suppose that $(t_0,x_0)
\in {\cal G}$. An application of the implicit function theorem to
$F(s,x,t) := \xi(t|s,x) - g(t)$ at the point
$(t_0,x_0,\xi(f(t_0,x_0)))$ implies that $f(s,x) < \infty$ in some
neighborhood $U$ of $(t_0,x_0)$ and that $f\big|_U \in C^2(U,\R)$.
The continuity of $g$ implies that $\{(s,x): x < g(t)\}$ is a
neighborhood of $(t_0,x_0)$, and the fact that $g$ is $C^2$ and
$f$ is continuous on $U$ implies that $\{(s,x):
g'(f(s,x))-I+\gamma g(f(s,x)) < 0 \}$ is a neighborhood of
$(t_0,x_0)$. \hfill $\Box$

\s

\n {\bf Proof of Proposition \ref{deviations}.} The compactness of $G$ implies the existence of
$\delta_1 > 0$ and $\delta_2 \in(0,\delta]$ such that
     $$
     \xi(t|t_0,x_0) < g(t) - \delta_1 \quad \mbox{ for } t_0 \le t \le f(t_0,x_0) - \delta
     $$
and
     $$
     \xi(t|t_0,x_0) > g(t) + \delta_1 \quad \mbox{ for } t  =
     f(t_0,x_0)+\delta_2
     $$
for all $(t_0,x_0) \in G$.  These two inequalities imply that
\begin{align*}
  \PP^{t_0,x_0}(|\tau^\e - f(t_0,x_0)| &> \delta) \\
    & \le  \PP^{t_0,x_0}\left(\{\tau^\e < f(t_0,x_0)-\delta\} \cup \{\tau^\e > f(t_0,x_0)+\delta_2\}\right)\\
    & \le   \PP^{t_0,x_0}\left(\{X_s^\e \ge g(s) \mbox{ for some }s \in [t_0,f(t_0,x_0)-\delta)\}\right. \\
      &  \hspace{8ex} \cup \left.\{X_{(f(t_0,x_0)+\delta_2)}^\e < g(f(t_0,x_0)+\delta_2)\}\right)\\
    & \le   \PP^{x_0,t_0}\left(\sup_{t_0 \le s \le f(t_0,x_0)+\delta_2}
     |X_s^\e - \xi(s|t_0,x_0)| \ge \delta_1\right).
\end{align*}
From (\ref{diff}) and the definition of $\xi(s|t_0,x_0)$ we have
 $$
|X_s^\e - \xi(s|t_0,x_0)| \leq  \gamma \int_{t_0}^s|X^\e(u) -
\xi(u|t_0,x_0)| du + \e |W_s-W_{t_0}|,
  $$
and then Gronwall's inequality gives
   $$
\sup_{t_0 \le s \le t} |X_s^\e - \xi(s|t_0,x_0)| \leq \e
e^{\gamma(t-t_0)} \sup_{t_0 \le s \le t}|W_s-W_{t_0}|.
  $$
Define $T = \sup\{f(t_0,x_0)+\delta_2 - t_0: (t_0,x_0) \in G\}$, and note that the
compactness of $G$ implies that $T < \infty$.  We have
\begin{eqnarray*}
\sup_{t_0 \le s \le f(t_0,x_0)+\delta_2} |X_s^\e -
\xi(s|t_0,x_0)|
   &\le& \e e^{\gamma T} \sup_{t_0 \le s \le
t_0+T}|W_s-W_{t_0}| \\
& \stackrel{\mbox{dist}}{=}& \e e^{\gamma T}
\sup_{0 \le s \le T}|W_s|.
\end{eqnarray*}
Therefore
 \begin{eqnarray*}
  \PP^{t_0,x_0}(|\tau^\e - f(t_0,x_0)| > \delta)
      & \le & \PP^{x_0,t_0}\left(\sup_{t_0 \le s \le f(t_0,x_0)+\delta_2}
     |X_s^\e - \xi(s|t_0,x_0)| \ge \delta_1 \right)\\
     & \le & \PP\left(\e e^{\gamma T}\sup_{0 \le s \le T}|W_s| \ge \delta_1\right)\\
     & \le & 2 \PP\left( |W_T| \ge \frac{\delta_1 e^{-\gamma T}}{\e}\right)
     \end{eqnarray*}
and the result follows directly.  \hfill $\Box$

\s

The proof of Theorem \ref{thm-unif} is rather lengthy so we will split the remainder of this section up into several subsections highlighting the main components which will be tied together in Section \ref{sec-proof}. We begin with a simple example for motivational purposes.

\subsection{Non-leaky case with constant threshold} \label{sec-special}

Suppose that $\gamma =0$ and $I > 0$ and that the threshold
$g(t)= B$ is constant. Then $X_t^\e$ is just Brownian motion with
drift $I$ and ${\cal G} = \{(s,x): x < B\}$ with $f(t) =
t_0+(B-x)/I$.  It can be shown (see for instance \cite{bnr}) that
for any $x_0 < B$ the first-passage time density is given
explicitly by
\begin{equation*}
p^\e(t|t_0,x_0) = \frac{B-x_0}{t-t_0} q^\e(t|t_0,x_0)
\end{equation*}
for $t \geq t_0$ where $q^\e(t|t_0,x_0)$ is the transition density
at $B$ for the random variable $X_t^\e$ given $X_{t_0}^\e = x_0$.
Thus
\begin{eqnarray*}
p^\e(t|t_0,x_0) & = &
\frac{B-x_0}{\sqrt{2\pi}\e(t-t_0)^{3/2}}\exp\left\{-\frac{(B - x_0
- I(t-t_0))^2}{2\e^2(t-t_0)}\right\}   \\
  & = & \frac{I(f(t_0,x_0)-t_0))}{\sqrt{2\pi}\e(t-t_0)^{3/2}}\exp\left\{-\frac{I^2(t-f(t_0,x_0))^2}{2\e^2(t-t_0)}\right\}.
\end{eqnarray*}
Replacing $t$ by $f(t_0,x_0) + \e t$ we get
\begin{align*}
 \e p^\e(\e t+ f(t_0,x_0)|t_0,x_0)
  &=
\frac{I(f(t_0,x_0)-t_0))}{\sqrt{2\pi}(f(t_0,x_0)-t_0+ \e t)^{3/2}} \\
& \quad \times \exp\left\{-\frac{I^2 t^2}{2(f(t_0,x_0)-t_0+ \e t)}\right\}
\\[2ex]
 & \to
   \frac{I(f(t_0,x_0)-t_0))}{\sqrt{2\pi}(f(t_0,x_0)-t_0)^{3/2}}
\exp\left\{-\frac{I^2 t^2}{2(f(t_0,x_0)-t_0)}\right\} \\[2ex]
  & =  p_\tau(t|t_0,x_0)
\end{align*}
as $\e \to 0$ because in this setting $\sigma^2(t_0,x_0) =
f(t_0,x_0)-t_0$ and $m(t_0,x_0) = I$.  This shows the pointwise
convergence implied by (\ref{clt}).  The full strength of
(\ref{clt}) will follow from the techniques developed below for
the general case.

\subsection{Durbin's Theorem} \label{sec-general}

If $\gamma > 0$ or $g$ is not constant, then we no longer have
explicit formulas at our disposal. However, since $X_t^\e$ is a
Gaussian process, we have the following result of Durbin
\cite[page 100]{jd}, valid for any $\gamma \ge 0$ and any $C^2$
function $g$.

\begin{thm} \label{thm-durbin} Suppose that $x_0 < g(t_0)$.  For $t > t_0$
\begin{equation} \label{peps}
p^\e(t|t_0,x_0) = b^\e(t|t_0,x_0)q^\e(t|t_0,x_0)
\end{equation}
where $q^\e(t|t_0,x_0)$ is the density under $\PP^{t_0,x_0}$ of
$X_t^\e$ evaluated at $g(t)$, and
\begin{equation*}
 b^\e(t|t_0,x_0):= \lim_{s\nearrow t}\frac{1}{t-s}\mathbb{E}^{t_0,x_0}[\mathbf{1}_{\tau^\e>s}(g(s)-X_s^\e)|X_t^\e=g(t)].
\end{equation*}
 \end{thm}

We call $q^\e$ the \emph{density} term and $b^\e$ the \emph{slope}
term in the decomposition \eqref{peps}.

\subsection{Analysis of the density term} \label{sec-dens}

We deal with $q^\e$ in much the same way as we dealt with $p^\e$
in the example from Section \ref{sec-special}.  Define
$$
\sigma^2(t|t_0) = \int_{t_0}^t e^{-2\gamma(t-s)}\,ds =
\left\{\begin{array}{cl}(1-e^{-2\gamma(t-t_0)})/2\gamma &
\hspace{.3in} \mbox{ if }\gamma > 0, \\[0.5ex]
t-t_0 & \hspace{.3in}
\mbox{ if }\gamma = 0. \end{array}\right.
   $$
Then
 \begin{equation} \label{oudens1}
   q^\e(t|t_0,x_0) = \frac{1}{\sqrt{2 \pi}\e \sigma(t|t_0)}
\exp\{-\frac{(g(t) - \xi(t|t_0,x_0))^2}{2\e^2\sigma^2(t|t_0)}\}
\end{equation}
where $\xi(t|t_0,x_0)$ is the solution to the ODE $x' = -\gamma x
+ I$ for $t \ge t_0$ with $\xi(t_0|t_0,x_0) = x_0$. Note that $\sigma_\tau^2(t_0,x_0)= \sigma^2(f(t_0,x_0)|t_0)/m^2(t_0,x_0)$.

\begin{lemma} \label{oudenslem}
If $G$ is a compact subset of ${\cal G}$ then there exist $\delta,
\e_0,K,\sigma_1 > 0$ so that
\begin{equation} \label{mqpest}
  \big|m(t_0,x_0)\e q^\e( f(t_0,x_0)+u|t_0,x_0) -  p_\tau(u/\e|t_0,x_0)\big| \leq K \e e^{-u^2/2 \e^2\sigma_1^2}
\end{equation}
and
  \begin{equation} \label{qest}
  (\e+|u|) q^\e(f(t_0,x_0)+u |t_0,x_0) \leq K  e^{-u^2/2\e^2\sigma_1^2}
\end{equation}
for all $(t_0,x_0) \in G$, $|u| \le \delta$ and $\e < \e_0$.
\end{lemma}

\s

\n{\bf Proof.}  The compactness of $G$ implies the existence of
$\delta_1
> 0$ such that $f(t_0,x_0) \ge t_0 +2 \delta_1$ for all $(t_0,x_0)
\in G$.  Noting that
$$g(f(t_0,x_0))-\xi(f(t_0,x_0)|t_0,x_0) = 0$$
and
$$g'(f(t_0,x_0))-\xi'(f(t_0,x_0)|t_0,x_0) = -m(t_0,x_0)$$
we obtain by Taylor's theorem
  $$
  \frac{\left[g(f(t_0,x_0)+u)-\xi(f(t_0,x_0)+u|t_0,x_0)\right]^2}{\sigma^2(f(t_0,x_0)+u |t_0)}
   = \frac{u^2}{\sigma_\tau^2(t_0,x_0)}\left(1+uR_1(t_0,x_0,u)\right)
   $$
where the remainder term satisfies
  $$
  |R_1(t_0,x_0,u)| \le K_1 \quad \quad \mbox{ for all }(t_0,x_0) \in G \mbox{ and }|u| \le \delta_1
  $$
for some $K_1$.  Similarly we have
  $$
  \frac{m(t_0,x_0)}{\sigma(f(t_0,x_0)+u |t_0)}
   = \frac{1}{\sigma_\tau(t_0,x_0)}\left(1+uR_2(t_0,x_0,u)\right)
   $$
where the remainder term satisfies
  $$
  |R_2(t_0,x_0,u)| \le K_2 \quad \quad \mbox{ for all }(t_0,x_0) \in G \mbox{ and }|u| \le \delta_1
  $$
for some $K_2$.  For ease of notation in the following calculation
we drop the arguments of $\sigma_\tau$ and $R_1$ and $R_2$.  We
have
%\hspace{1em}
  \begin{eqnarray*}
  \lefteqn{
  \left|m(t_0,x_0)\e q^\e(f(t_0,x_0)+u |t_0,x_0)-
  p_\tau(u/\e|t_0,x_0)\right|} \\[2ex]
  & = &
  \frac{1}{\sqrt{2
  \pi}\sigma_\tau}\left|(1+uR_2)e^{-u^2(1+uR_1)/2\e^2
  \sigma_\tau^2} - e^{-u^2/2\e^2
  \sigma_\tau^2}\right|\\[2ex]
  & \le &
  \frac{1}{\sqrt{2
  \pi}\sigma_\tau}\left( |uR_2|e^{-u^2(1+uR_1)/2\e^2
  \sigma_\tau^2} + \left|e^{-u^2(1+uR_1)/2\e^2
  \sigma_\tau^2} - e^{-u^2/2\e^2
  \sigma_\tau^2}\right|\right)\\[2ex]
    & \le &
  \frac{1}{\sqrt{2
  \pi}\sigma_\tau}\left(|uR_2|e^{-u^2(1+uR_1)/2\e^2
  \sigma_\tau^2}  \right. \\
  & & \hspace{1.2in} + \left. \frac{|u^3R_1|}{2 \e^2 \sigma_\tau^2} \max\left\{e^{-u^2(1+uR_1)/2\e^2
  \sigma_\tau^2}, e^{-u^2/2\e^2
  \sigma_\tau^2} \right\}\right).\end{eqnarray*}
At this point choose $\delta \le \delta_1$ so that $K_1 \delta \le
1/3$.  Then for $|u| \le \delta$ we have
   \begin{eqnarray*}
  \lefteqn{
  \left|m(t_0,x_0)\e q^\e(f(t_0,x_0)+u |t_0,x_0)-
  p_\tau(u/\e|t_0,x_0)\right|}  \\[2ex]
    & \le &
  \frac{1}{\sqrt{2
  \pi}\sigma_\tau}\left( K_2|u|e^{-u^2/3\e^2
  \sigma_\tau^2} + \frac{K_1|u^3|}{2 \e^2 \sigma_\tau^2} e^{-u^2/3\e^2
  \sigma_\tau^2} \right)\\[2ex]
  &= & \frac{\e}{\sqrt{2
  \pi}}\left( K_2\left|\frac{u}{\e \sigma_\tau}\right|  + \frac{K_1}{2}\left|\frac{u}{\e \sigma_\tau}\right|^3\right) e^{-u^2/3\e^2
  \sigma_\tau^2}\\[2ex]
  & \le & \frac{\e}{\sqrt{2
  \pi}}\left( K_2  + \frac{K_1}{2}\right)K_3 e^{-u^2/4\e^2
  \sigma_\tau^2}\end{eqnarray*}
where
  $$
    K_3 = \max\left( \max\{|x| e^{-x^2/12}: x\in \R\}, \max\{|x|^3
    e^{-x^2/12}: x \in\R\} \right).
  $$
The result (\ref{mqpest}) follows directly, with $\sigma_1^2 =
2\max\{\sigma_\tau^2(t_0,x_0): (t_0,x_0) \in G\}$.  The proof of
(\ref{qest}) is similar (and simpler) and is left to the reader.
\hfill $\Box$

\subsection{Analysis of the slope term} \label{sec-slope}

The slope term
   \begin{equation*}
 b^\e(t|t_0,x_0):= \lim_{s\nearrow t}\frac{1}{t-s}\mathbb{E}^{t_0,x_0}[\mathbf{1}_{\tau^\e>s}(g(s)-X_s^\e)|X_t^\e=g(t)]
\end{equation*}
involves the pinned process $X^\e$ given $X^\e_{t_0} = x_0$ and
$X^\e_t = y$.  For $t_0 < t$ define
\begin{equation} \label{psi}
  \psi(s|t_0,t)  = \left\{
    \begin{array}{ll}
      \displaystyle{\frac{\sinh \gamma(s-t_0)}{\sinh \gamma(t-t_0)}} & \mbox{ if } \gamma >
      0
      \\[2ex]
      \displaystyle{\frac{s-t_0}{t-t_0}} & \mbox{ if } \gamma = 0
      \end{array} \right.
     \end{equation}
and\begin{equation} \label{m}
  \mu(s|t_0,x_0,t,y) = \xi(s|t_0,x_0) + \psi(s|t_0,t)(y -
\xi(t|t_0,x_0))
  \end{equation}
   for $t_0 \le s \le t$. Our next proposition gives us a useful representation for the pinned process.

\begin{proposition} \label{prop-pinned}  For $t_0 < t$ the
conditional distribution of $\{X^\e_s: t_0 \le s \le t\}$ given
$X^\e_{t_0} = x_0$ and $X^\e_t = y$ can be written
  $$
  X^\e_s = \mu(s|t_0,x_0,t,y) + \e U_s, \quad \quad s \in [t_0,t]
  $$
where $\E (U_s) = 0$ and
  $$\mathbb{E}\left(\sup_{t_0 \le s \le t}|U_s|\right) \le
     \left\{ \begin{array}{cc}
       K \sqrt{(e^{2\gamma(t- t_0)}-1)/2 \gamma} & \mbox{
       if }\gamma > 0 \\[1ex]
       K (t-t_0) & \mbox{ if } \gamma = 0.
       \end{array} \right.
  $$
\end{proposition}

\s

The proof of Proposition \ref{prop-pinned} relies on the following
standard result regarding the conditioned law of a Gaussian
process.  We include the proof for completeness.

\begin{lemma} \label{lem-gauss} Suppose that $\{Z_s:  s \in S\}$ is
a real valued Gaussian process defined on some set $S \subset \R$
with mean $\mu(s)$ and covariance $\rho(s,t)$. Let $ N \ge 1$ and
$t_1, t_2, \ldots t_N \in S$, and suppose that the matrix $A :=
(\rho(t_i,t_j))_{i,j=1,...,N}$ is invertible with $A^{-1} = B$.
Then the conditional law of $\{Z_s : s \in S\}$ given
$Z_{t_1}=z_1,...,Z_{t_N}=z_N$ is given by
   $$
\mu(s) + \sum_{i,j=1}^N  \rho(s,t_i)B_{ij}(z_j-\mu(t_j)) + U_s,
\quad \quad s \in S,
$$
where
$$
U_s = (Z_s - \mu(s)) - \sum_{i,j=1}^N
\rho(s,t_i)B_{ij}(Z_{t_j}-\mu(t_j)).
$$
The process $\{U_t: s \in S \}$ is Gaussian with mean 0 and
covariance function $\tilde{\rho}(s,t) = \rho(s,t) -
\sum_{i,j=1}^N \rho(s,t_i)B_{ij}\rho(t,t_j)$.
\end{lemma}

\n {\bf Proof.}  Define the process
$$
V_s = \mu(s) + \sum_{i,j=1}^N  \rho(s,t_i)B_{ij}(Z_{t_j}-\mu(t_j))
$$
for $s \in S$.  Clearly $\{V_s: s \in S\}$ is measurable with
respect to $\sigma\{Z_{t_1},..,Z_{t_N}\}$ and is a Gaussian
process.  Moreover $V_{t_j} = Z_{t_j}$ for $j=1,2,...,N$. Now
define $U_s = Z_s - V_s$.  The process $\{U_s: s \in S\}$ is a
mean zero Gaussian process, and a direct algebraic calculation
gives that $\mbox{Cov}(U_s,Z_{t_j}) =0$ for all $j = 1, 2,
\ldots,N$, so that $\{U_s: s \in S\}$ is independent of
$\sigma\{Z_{t_1},..,Z_{t_N}\}$.  Therefore the law of $Z_s = U_s +
V_s$ given $Z_{t_i} = z_i$, for $i=1,...,N$ is the same as the law
of $U_s + \mu(s) + \sum_{i,j=1}^N \rho(s,t_i)B_{ij}(z_j-\mu(t_t) $, and the first
assertion is proved. The covariance of $\{U_s: s \in S\}$ is a
direct algebraic calculation and is left to the reader.
   \hfill $\Box$

\s

\n {\bf Proof of Proposition \ref{prop-pinned}.}  We will apply
Lemma \ref{lem-gauss} to the process $\{X_s^\e: s \ge t_0\}$ given by
   $$
   dX_s^\e = (-\gamma X_s^\e +I)\,ds + \e dW_s
   $$
with initial condition $X_{t_0}^\e = x_0$, conditioned by $X_t^\e = y$.
We will take $N = 1$ and condition at the single point $t$, so that
$N = 1$ and $B_{11} = \rho(t,t)^{-1}$.  For $s \ge t_0$ we have
$\mu(s) = \xi(s|t_0,x_0)$, and for $s_1,s_2 \ge t_0$ we have
    $$
    \rho(s_1,s_2) = \left\{
    \begin{array}{cl}
      (\e^2/2\gamma)\left(e^{-\gamma|s_1-s_2|} - e^{-\gamma(s_1+s_2-2
      t_0)}\right) & \mbox{ if } \gamma > 0 \\[1ex]
          \e^2 \min(s_1-t_0,s_2-t_0) & \mbox{ if } \gamma = 0.
            \end{array} \right.
            $$
Lemma \ref{lem-gauss} implies that the conditional law of $\{X_s^\e:
t_0 \le s \le t\}$ given $X_t^\e = y$ is
   $$
   \mu(s|t_0,x_0,t,y) + U^\e_s
   $$
where
  $$
  \mu(s|t_0,x_0,t,y) =  \xi(s|t_0,x_0) + \frac{\rho(s,t)}{\rho(t,t)}(y-\xi(t|x_0,t_0)
       $$
and
$$
U^\e_s = X_s^\e - \xi(s|t_0,x_0) -
\frac{\rho(s,t)}{\rho(t,t)}(X_t^\e-\xi(t|t_0,x_0)).
  $$
\n {\bf Case 1:} $\gamma = 0$.  We have $X_s^\e - \xi(s|t_0,x_0) =
\e(W_s-W_{t_0})$ and $\rho(s_1,s_2) = \e^2 \min(s_1-t_0,s_2-t_0)$,
so that $U^\e_s = \e U_s$ where
   $$
U_s = \left[W_s-W_{t_0} -
\left(\frac{s-t_0}{t-t_0}\right)(W_t-W_{t_0})\right].
  $$
The process $\{U_s: t_0 \le s \le t\}$ is the Brownian bridge on the
interval $[t_0,t]$.  For $0 \le u \le 1$ define $\widetilde{W}_u =
(W_{t_0+u(t-t_0)}-W_{t_0})/\sqrt{t-t_0}$, then $\widetilde{W}$ is
standard Brownian motion and
  $$
U_s = \sqrt{t-t_0}\left[\widetilde{W}_{(s-t_0)/(t-t_0)} -
\left(\frac{s-t_0}{t-t_0}\right)\widetilde{W}_1 \right].
  $$
We have
  $$
  \sup_{t_0 \le s \le t}|U(s)|
   =  \sqrt{t-t_0} \sup_{0 \le u \le 1}|\widetilde{W}_u-u\widetilde{W}_1|
   $$
and in particular
  $$
  E\left(\sup_{t_0 \le s \le t}|U(s)|\right)
   =  \sqrt{t-t_0} E\left(\sup_{0 \le u \le 1}|\widetilde{W}_u-u\widetilde{W}_1|\right)= K \sqrt{t-t_0}.
   $$

\s

\n {\bf Case 2:} $\gamma >0$.  For $t_0 \le s \le t$ we have
   $$
   X_s -\xi(s|t_0,x_0) = \e \int_{t_0}^s e^{-\gamma(s-v)}dW_v
     $$
and
     $$ \frac{\rho(s,t)}{\rho(t,t)}
     = \frac{\left(e^{-\gamma(t-s)} - e^{-\gamma(s+t-2t_0)}\right)}{\left(1 - e^{-2\gamma(t-t_0)}\right)}
     %=\frac{\sinh \gamma(s-t_0)}{\sinh \gamma(t-t_0)}
     .
      $$
Therefore $U^\e_s = \e U_s$ where
  $$
U_s    =   \int_{t_0}^s e^{-\gamma(s-v)}dW_v -
\frac{\left(e^{-\gamma(t-s)} -
e^{-\gamma(s+t-2t_0)}\right)}{\left(1 -
e^{-2\gamma(t-t_0)}\right)} \int_{t_0}^t e^{-\gamma(t-v)}dW_v.
   $$
and so
  $$
e^{\gamma s} U_s    =   \int_{t_0}^s e^{\gamma v}dW_v -
\frac{\left(e^{-2\gamma(t-s)} - e^{-2
\gamma(t-t_0)}\right)}{\left(1 - e^{-2\gamma(t-t_0)}\right)}
\int_{t_0}^t e^{\gamma v}dW_v.
   $$
For fixed $t_0 < t$ define $a:[t_0,t] \to \R$ by
   $$
   a(s) = \frac{\left(e^{-2\gamma(t-s)} - e^{-2
\gamma(t-t_0)}\right)}{\left(1 - e^{-2\gamma(t-t_0)}\right)}.
$$
The function $a$ is continuous and strictly increasing from
$[t_0,t]$ onto $[0,1]$, with inverse
  $$
  a^{-1}(u) = t_0 + \frac{1}{2\gamma} \log\left((1-u)+u e^{2 \gamma (t-t_0)}\right), \quad \quad 0 \le u \le 1.
  $$
For $0 \le u \le 1$ define
   $$
   \widehat{W}(u) = \sqrt{\frac{2 \gamma}{e^{2 \gamma t}-e^{2 \gamma t_0}}} \cdot  \int_{t_0}^{a^{-1}(u)} e^{\gamma
   v }dW_v
   $$
Then $\{\widehat{W}(u): 0 \le u \le 1\}$ is a standard Brownian
motion, and
 $$
e^{\gamma s} U_s = \sqrt{\frac{e^{2 \gamma t}-e^{2 \gamma t_0}}{2
\gamma}} \left( \widehat{W}(a(s)) - a(s)\widehat{W}(1)\right).
   $$
Therefore
    \begin{eqnarray*}
    \sup_{t_0 \le s \le t}|U(s)|
      &  \le &  e^{-\gamma t_0}\sup_{t_0 \le s \le t}|e^{\gamma
      s}U(s)|\\
      & = & \sqrt{\frac{e^{2\gamma(t- t_0)}-1}{2 \gamma}} \sup_{t_0 \le s \le t}\left|\widehat{W}(a(s))  - a(s) \widehat{W}(1)\right|\\
   & = & \sqrt{\frac{e^{2\gamma(t- t_0)}-1}{2 \gamma}} \sup_{0 \le u \le 1}\left|\widehat{W}(u)  - u
   \widehat{W}(1)\right|.
    \end{eqnarray*}
In particular we have
  $$
  E\left(\sup_{t_0 \le s \le t}|U(s)|\right)
   \le K \e \sqrt{\frac{e^{2\gamma(t- t_0)}-1}{2 \gamma}}.
   $$  \eopt

\s

\s

For $t_0 < t$ define
    \begin{equation} \label{beta}
  \beta(s|t_0,x_0,t) = g(s)-\mu(s|t_0,x_0,t,g(t)) \quad \quad t_0
  \le s \le t,
  \end{equation}
so that $\beta(s| t_0,x_0,t,)$ is the gap between the threshold
$g(s)$ and the conditional expected value of $X^\e_s$ given
$X_{t_0}^\e = x_0$ and $X_t^\e = g(t)$. A direct calculation gives
    \begin{eqnarray}
  \beta'(s|t_0,x_0,t)
    & = & g'(s)+\gamma \xi(s|t_0,x_0)-I - \psi'(s|t_0,t)(g(t)-\xi(t|t_0,x_0)  \label{betaderiv}\\[2ex]
    & = & \left(g'(s)-\frac{g(t)-x_0}{t-t_0}\right) +\left(\frac{g(t)-x_0}{t-t_0}\right)
     \left(1-\frac{\gamma (t-t_0) \cosh \gamma(t-s)}{\sinh \gamma(t-t_0)}\right) \nonumber \\
  & & \mbox{} + ( \gamma g(t)-I)\left(\frac{ \cosh \gamma(t-s) - \cosh \gamma (s-t_0)}{\sinh
  \gamma(t-t_0)}\right)\label{betaderiv2}
  \end{eqnarray}
if $\gamma > 0$ and
    \begin{equation}
  \beta'(s|t_0,x_0,t)
    = \left(g'(s)-\frac{g(t)-x_0}{t-t_0}\right) \label{betaderiv3}
  \end{equation}
if $\gamma = 0$.

\s

Following Durbin \cite{jd} we can write
\begin{equation} \label{bdecomp}
b^\e(t|t_0,x_0) = b_1^\e(t|t_0,x_0) - \bar{b}^\e(t|t_0,x_0)
\end{equation}
where
$$b_1^\e(t|t_0,x_0) := \lim_{s\nearrow
t}\frac{1}{t-s}\mathbb{E}^{t_0,x_0}[(g(s)-X_s^\e)|X_t^\e=g(t)]$$
and
$$\bar{b}^\e(t|t_0,x_0) := \lim_{s\nearrow t}(t-s)^{-1}\mathbb{E}^{t_0,x_0}[\mathbf{1}_{\tau^\e\leq s}(g(s)-X_s^\e)|X_t^\e=g(t)].
 $$
By Proposition \ref{prop-pinned} we have
\begin{eqnarray} \label{b1exp}
b_1^\e(t|t_0,x_0)
 &=& \lim_{s\nearrow t}\frac{\beta(s|t_0,x_0)}{t-s} \nonumber \\
 &=& -\beta'(t|t_0,x_0,t) \nonumber \\
&=& -g'(t)-\gamma \xi(t|t_0,x_0) + I -
\psi'(t|t_0,t)(g(t)-\xi(t|t_0,x_0)).
 \end{eqnarray}
In particular $b_1^\e(t|t_0,x_0)$ is independent of $\e$, and
henceforth we shall write $b_1(t|t_0,x_0) = b_1^\e(t|t_0,x_0)$.
Notice that
  \begin{equation} \label{b1f}
   b_1(f(t_0,x_0)|t_0,x_0) = -g'(f(t_0,x_0)) -\gamma g(f(t_0,x_0)) + I =
   m(t_0,x_0).
   \end{equation}

\s

The following result of Durbin \cite{jd} will enable us to control
the $\bar{b}^\e$ term in \eqref{bdecomp}.

\begin{proposition} \label{prop-bbar}
For any $t > t_0$ we have
\begin{eqnarray}
\bar{b}^\e(t|t_0,x_0) &=& \int_{t_0}^t
b_1(t|r,g(r))\bar{p}^\e(r|t_0,x_0,t,g(t))dr. \label{bbarint}
\end{eqnarray}
where $\bar{p}^\e(r|t_0,x_0,t,y)$ denotes the conditional density
function for $\tau^\e$ given that $X_t^\e=y$ and $X_{t_0}^\e =
x_0$
\end{proposition}

\n {\bf Proof:} This a corrected restatement of \cite[Section 6,
equation (27)]{jd}.  The Markov property for $X^\e$ gives
      \begin{eqnarray*}
   \bar{b}^\e(t|t_0,x_0) &=& \lim_{s \nearrow t} \int_{t_0}^s
 \E^{r,g(r)}\left(\left.  \frac{g(s)-X^\e(s)}{t-s} \right| X^\e(t) = g(t) \right) \bar{p}^\e(r|t_0,x_0,t,g(t))dr \\[2ex]
     & = & \lim_{s \nearrow t} \int_{t_0}^s
 \left(  \frac{ \beta(s|r,g(r),t)}{t-s} \right)
 \bar{p}^\e(r|t_0,x_0,t,g(t))dr.
  \end{eqnarray*}
It is easily checked from the expressions (\ref{betaderiv2}) and
(\ref{betaderiv3}) that $|\beta'(s|r,g(r),t)|$ is bounded for $t_0
\le r < s < t$, and since we also have $\beta(t|r,g(r),t)=0$, the passage of $\lim_{s \nearrow t}$ inside
the integral is justified by the bounded convergence theorem.
\eopt

\begin{remark} \label{rem-bnr-durbin}
{\rm   The substitution of (\ref{bbarint}) into (\ref{bdecomp})
gives an integral expression for $b^\e(t|t_0,x_0)$ in terms of the
conditional density function $\bar{p}^\e(r|t_0,x_0,t,y)$.  The
relation (6.1) then gives an integral equation for the first passage
density $p^\e(t|t_0,x_0)$ which is a special case of the integral
equation of Buonocore, Nobile and Ricciardi \cite{bnr}. In
particular \cite{bnr} can be used to give an alternative derivation
of (\ref{b1exp}).  The paper \cite{bnr} deals with time-homogeneous
diffusion processes, and the paper \cite{jd} deals with Gaussian
processes, and we are working in the intersection of these two
classes of processes.}
\end{remark}

\begin{proposition} \label{prop-blim} Let $G$ be a compact
subset of ${\cal G}$. Then there exist $\delta
> 0$ and $K < \infty$ such that
  \begin{equation} \label{blim}
  \left|b^\e(t|t_0,x_0) - b_1(t|t_0,x_0)\right| \le \e K
  \end{equation}
whenever $(t_0,x_0) \in G$ and $|t-f(t_0,x_0)| \le  \delta$.
\end{proposition}

We prepare for the proof with the following pair of Lemmas.

\begin{lemma}  \label{lem-gmdiff} Given the compact set $G \subset {\cal G}$ there are positive $\delta$, $\delta_1$, $k_1$ and $k_2$ such that
  $$
  \beta(s|t_0,x_0,t)  \ge
    \left\{\begin{array}{cl} k_1(t-s) & \quad \mbox{ if } t-\delta_1 \le s \le t
    \\[1ex]
       k_2 & \quad \mbox{ if } t_0 \le s \le t-  \delta_1
    \end{array} \right.
        $$
whenever $(t_0,x_0) \in G$ and $|t-f(t_0,x_0)| \le \delta$.
   \end{lemma}

\n{\bf Proof.}  We give the proof for the case $\gamma > 0$; the
proof for the case $\gamma =0$ is essentially the same. Equation
(\ref{betaderiv2}) shows that $\beta'(s|t_0,x_0,t)$ is a
continuous function of $(t_0,x_0,t,s) $ on the set where $t >
t_0$.  Also, putting $s = t = f(t_0,x_0)$ in (\ref{betaderiv})
gives
  $$
  \beta'(f(t_0,x_0)|t_0,x_0,f(t_0,x_0)) = g'(f(t_0,x_0))+\gamma g(f(t_0,x_0)) - I =  - m(t_0,x_0) <0.
  $$
for $(t_0,x_0) \in {\cal G}$.  The compactness of $G$ gives
$\delta_2> 0$ and $k_1 > 0$ such that $\beta'(s|t_0,x_0,t)\le -k_1
$ whenever $(t_0,x_0) \in G$ and $f(t_0,x_0)-\delta_2  \le s \le t \le
f(t_0,x_0)+ \delta_2$.  Since $\beta(t|t_0,x_0,t) = 0$ it follows
that
 \begin{equation}
 \label{beta1}\beta(s|t_0,x_0,t) \ge k_1(t-s)
 \end{equation}
whenever $(t_0,x_0) \in G$ and $f(t_0,x_0) - \delta_2 \le s \le t
\le f(t_0,x_0)+\delta_2$.

The definition of $f(t_0,x_0)$ as the time of first intersection
of $\xi(t|t_0,x_0)$ with $g(t)$ implies that
  $$\min\{\beta(s|t_0,x_0,f(t_0,x_0)): t_0 \le s \le
  f(t_0,x_0)-\delta_2 \}> 0
  $$
for all $(t_0,x_0) \in G$.  The compactness of $G$ gives the
existence of $\delta_3 > 0$ and $k_2 > 0$ such that
  \begin{equation} \label{beta2}
  \beta(s|t_0,x_0,t) \ge k_2
  \end{equation}
whenever $(t_0,x_0) \in G$ and $|t-f(t_0,x_0)| \le \delta_3$ and
$t_0 \le s \le f(t_0,x_0)-\delta_2$.  The result is now a simple
consequence of (\ref{beta1}) and (\ref{beta2}), with $\delta =
\min(\delta_2/2,\delta_3)$, $\delta_1 = \delta_2/2$. \hfill $\Box$

\begin{lemma} \label{lem-b1} Let $H$ be a compact subset of $\R$.  There is $K_1$ such that
   $$
   |b_1(t|r,g(r))| \le K_1(t-r)
   $$
whenever $r, t \in H$ and $r < t$.
 \end{lemma}

\n{\bf Proof.}  Again we give the proof for the case $\gamma
> 0$ and leave $\gamma = 0$ to the reader.  Putting $s = t$ and $(t_0,x_0) = (r,g(r))$ in (\ref{betaderiv2})
gives
\begin{align*}
b_1(t|r,g(r))
&= \left(\frac{g(t)-g(r)}{t-r}  -g'(t)\right) + \left(\frac{\gamma(t-r)}{\sinh[\gamma(t-r)]}-1\right)\left(\frac{g(t)-g(r)}{t-r}\right) \\
& \quad + \gamma
\left(\frac{\cosh[\gamma(t-r)]-1}{\sinh[\gamma(t-r)]}\right)(g(t)-I).
\end{align*}
Using the inequalities
  $$
  \left|\frac{\cosh u-1}{\sinh u}\right| \le \frac{|u|}{6}
  \quad \mbox{ and } \quad
   \left| \frac{u}{\sinh u} - 1\right| \le \min\left(\frac{u^2}{6},1\right)
   \le \frac{|u|}{\sqrt{6}}
  $$
we get
 \begin{eqnarray*}
\frac{|b_1(t|r,g(r))|}{t-r}
& \le & \left|\frac{g(t)-g(r) -(t-r) g'(t)}{(t-r)^2} \right| + \frac{\gamma}{\sqrt{6}}\left|\frac{g(t)-g(r)}{t-r}\right| \\
& & + \frac{\gamma^2}{2}|g(t)-I|\\
  & \le & \frac{1}{2} \sup_{r \le s \le t}|g''(s)| + \frac{\gamma}{\sqrt{6}}\sup_{r \le s \le t}|g'(s)| + \frac{\gamma^2}{2}|g(t)-I|,
\end{eqnarray*}
and the result follows from the compactness of $H$. \hfill $\Box$

\s

\n{\bf Proof of Proposition \ref{prop-blim}.}  Let $\delta$,
$\delta_1$, $k_1$ and $k_2$ be as in Lemma \ref{lem-gmdiff}, and
then $K_1$ as in Lemma \ref{lem-b1} with $H=\{t: t_0 \le t \le
f(t_0,x_0)+\delta \mbox{ for some } (t_0,x_0) \in G\}$.  Now fix
$(t_0,x_0) \in G$ and $t \in [f(t_0,x_0)-\delta,
f(t_0,x_0)+\delta]$ and apply Proposition \ref{prop-pinned} with
$y = g(t)$.  Under the conditions $X^\e_{t_0} = x_0$ and $X^\e_t =
g(t)$ we have
    $$
   \tau^\e = \inf\{s \ge t_0: \e U_s = \beta(s|t_0,x_0,t)\}.
   $$
Define $\|U\|:=\sup_{t_0 \le s \le t}|U_s|$. If $\e \|U\| = x <
\min(k_1\delta_1,k_2)$ then by Lemma \ref{lem-gmdiff} we have $\e
U_s < \beta(s|t_0,x_0,t)$ for all $s \le t-x/k_1$ and so $\tau^\e
\ge t-x/k_1$ and $|b_1(t|\tau^\e,g(\tau^\e))| \le K_1x/k_1= \e
K_1\|U\|/k_1$.  If $ \e \|U\| \ge \min(k_1\delta_1,k_2)$ we can
use the estimate $|b_1(t|\tau^\e,g(\tau^\e))| \le K_1(t-\tau^\e)
\le K(t-t_0)$ from Lemma \ref{lem-b1}. Together we have
 $$
 b_1(t|\tau^\e,g(\tau^\e)) \le
\e \|U\| \max\left(\frac{K_1}{k_1},
\frac{K_1(t-t_0)}{\min(k_1\delta_1,k_2)} \right).
 $$
Proposition \ref{prop-bbar} gives
     \begin{eqnarray*}
  |\overline{b}^\e(t|t_0,x_0)|
   &  \le & \int_{t_0}^t |b_1(t|r,g(r))|p^\e(r|t_0,x_0,t,g(t))\,dt\\
   & = & \mathbf{E}^{t_0,x_0}\big(|b_1(t|\tau^\e,g(\tau^\e))| \,\big| X^\e_t =
   g(t)\big)\\
   & \le & \e \mathbf{E}\|U\| \max\left(\frac{K_1}{k_1},
\frac{K_1(t-t_0)}{\min(k_1\delta_1,k_2)} \right),
   \end{eqnarray*}
and the result now follows by Proposition \ref{prop-pinned} and
the compactness of $G$. \hfill $\Box$

\subsection{Proof of Theorem \ref{thm-unif}} \label{sec-proof}

We are now ready to complete the proof of Theorem \ref{thm-unif}.
Given $G$, choose $\delta > 0$ sufficiently small and $K < \infty$
sufficiently large so that the results of Lemma \ref{oudenslem}
and Proposition \ref{prop-blim} are valid.  By Theorem
\ref{thm-durbin} we have
     \begin{eqnarray*}
 \lefteqn{\big|\e p^\e(f(t_0,x_0)+u|t_0,x_0) -
 p_\tau(u/\e|t_0,x_0) \big|} \\
   & \le & \big| m(t_0,x_0) \e q^\e(f(t_0,x_0)+u) -
    p_\tau(u/\e|t_0,x_0)\big|  \\
    & & \mbox{} + \big| b^\e(f(t_0,x_0)+u|t_0,x_0) -m(t_0,x_0)\big| \, \e q^\e(f(t_0,x_0)+u|t_0,x_0) \\
      &=& I + II.
  \end{eqnarray*}
Now by Proposition \ref{prop-blim}
\begin{eqnarray*}
\lefteqn{\big|b^\e(f(t_0,x_0)+u|t_0,x_0)-m(t_0,x_0)\big|} \\
   & \le & \big|b^\e(f(t_0,x_0)+u|t_0,x_0)-b_1(f(t_0,x_0)+u|t_0,x_0)\big|  \\
   & & \mbox{}+ \big|b_1(f(t_0,x_0)+u|t_0,x_0) - b_1(f(t_0,x_0)|t_0,x_0)\big|\\
   &  \le & \e K + |u|K_1
\end{eqnarray*}
where
   $$
   K_1 = \sup \{|b_1'(s|t_0,x_0)|: |s-f(t_0,x_0)| \le \delta \} <
   \infty.
   $$
(The finiteness of $K_1$ uses the fact that $f(t_0,x_0)-t_0$ is
bounded away from 0 on the compact set $G$.)  The result now
follows from Lemma \ref{oudenslem}, using (\ref{mqpest}) on $I$
and (\ref{qest}) on $II$.   \hfill $\Box$

\section{Proofs for Section \ref{sec-sifmtrans}} \label{sec-sifmtransproofs}

In what follows, $d$ denotes the standard quotient metric on
$\mathbb{S}$ induced by the Euclidean metric on $\mathbb{R}$. The
starting point for the proof of Theorem \ref{thm-outrans} is the
following splitting of $\mathbb{S}$.

\begin{proposition} \label{regions}
There exist neighborhoods $V_1 := B_{\delta_u}(\theta_u)$, $V_3 :=
B_{\delta_s}(\theta_s)$, and constants $\delta >0$, $N \in
\mathbb{N}$ such that
\begin{description}
\item[(1)]$d(\tilde{f}(\theta),V_1) > \delta$ for every $\theta
\notin V_1$

\item[(2)]$d(\tilde{f}(\theta),V_3^c) > \delta$ for every $\theta
\in V_3$

\item[(3)] For every $\theta \in V_2:=\mathbb{S}/(V_1 \cup V_3)$, we have
$\tilde{f}^n(\theta) \in V_3$, $\forall n \geq N$.
\end{description}
\end{proposition}

\s

\n {\bf Proof.} This is the same as Proposition 1 in
\cite{jm}. \eopt

\vs

We can write any $\phi \in B(\mathbb{S})$ in the form $\phi =
\phi_1 + \phi_2 + \phi_3$ where $\phi_i = \phi\mathbf{1}_{V_i}$.
The action of $T^\e$ can then be described by the block
decomposition
$$
T^\e = \left[ \begin{array}{ccc} T_{11}^\e & T_{12}^\e & T_{13}^\e \\
      T_{21}^\e & T_{22}^\e & T_{23}^\e \\T_{31}^\e & T_{32}^\e & T_{33}^\e \end{array} \right]
$$
where
\begin{equation} \label{Tij}
T_{ij}^\e \phi(\theta_0) =  \int \phi(\theta)
\tilde{p}_{ij}^\e(\theta|\theta_0) d\theta
\end{equation}
 with $\tilde{p}_{ij}^\e(\theta|\theta_0) =
\mathbf{1}_{V_i}(\theta_0)\mathbf{1}_{V_j}(\theta)\tilde{p}^\e(\theta|\theta_0)$.
The choice of the sets $V_1$, $V_2$ and $V_3$ implies that
   $$
T^0 = \left[ \begin{array}{ccc} T_{11}^0 & T_{12}^0 & T_{13}^0 \\
      0 & T_{22}^0 & T_{23}^0 \\0 & 0 & T_{33}^\e \end{array}
      \right].
$$
For $\e > 0$ write
 $$
T_{lp}^\e = \left[ \begin{array}{ccc} 0 &0 & 0\\
      T_{21}^\e & 0 & 0 \\T_{31}^\e & T_{32}^\e & 0 \end{array}
      \right]\quad \mbox{ and } \quad
       T_{up}^\e = \left[ \begin{array}{ccc} T_{11}^\e & T_{12}^\e & T_{13}^\e \\
      0 & T_{22}^\e & T_{23}^\e \\0 & 0 & T_{33}^\e \end{array}
      \right].
$$

\begin{lemma} \label{lower}  There are finite positive constants $K$ and $M$ such that $\|T^\e_{ij}\|_\infty \le
M\e e^{-K/\e^2}$ for $ij = 21$ or $31$ or $32$, and
$\|(T^\e_{22})^{N+1}\|_\infty \le  M\e e^{-K/\e^2}$.
\end{lemma}

\n {\bf Proof.}  For $ij = 21$ or $31$ or $32$ we have by
Proposition \ref{regions}
\begin{align*}
\|T_{ij}^\e \| = \sup_{\theta \in V_i} \mathbb{P}^\theta(\Theta_1^\e \in V_j)
      &\le\sup_{\theta \in V_i} \mathbb{P}^\theta(d(\Theta_1^\e, \tilde{f}(\theta)) > \delta) \\
     &\le \sup_{0 \le t \le 1}\mathbb{P}^t(|\tau^\e - f(t)| > \delta),
\end{align*}
and the first set of results follows by Proposition
\ref{deviations}. Also by Proposition \ref{regions}
  $$\|(T_{22}^\e)^{N+1} \|_\infty \le \sup_{\theta \in V_2} \mathbb{P}^\theta(d(\Theta_{N+1}^\e, \tilde{f}^{N+1}(\theta)) >
  \delta),
         $$
and the second result follows using Proposition \ref{deviations}
together with the inequality
  $$
  d(\Theta_{N+1}^\e,\tilde{f}^{N+1}(\theta)) \le \sum_{j=0}^N L^{N-i}d(\Theta^\e_{j+1}, \tilde{f}(\Theta^\e_j))
  $$
where $L = \sup |\tilde{f}'(\theta)|$.  \hfill $\Box$

\s

It follows directly from Lemma \ref{lower} that for any $r > 0$
there is $\e_0
> 0$ such that $\|T^\e_{lp}\|_\infty < r$ and all eigenvalues of
$T^\e_{22}$ have modulus less than $r$ whenever $\e < \e_0$.  In
order to complete the proof of Theorem \ref{thm-outrans} it
suffices to describe the eigenvalues of the operators $T_{11}^\e$
and $T_{33}^\e$ as $\e \to 0$.  This will be carried out in
Sections \ref{sec-local} and \ref{sec-unstable}.

\subsection{Behavior near a stable fixed point} \label{sec-local}

Recall that $T_{33}^\e$ is the restriction of the transition
operator $T^\e$ to a neighborhood $V_3$ of the stable fixed point
$\theta_s$, and that $f'(\theta_s) = c_s$.  The main result of
this section is

\begin{proposition} \label{prop-stabev}  Every non-zero eigenvalue
of $T_{33}^\e$ is of the form $c_s^n + O(\e)$ as $\e \to 0$ for
some $n \ge 0$.
\end{proposition}

We can reparameterize $\mathbb{S}$ so that $\theta_s = 0$ and $V_3
= (-\delta_s,\delta_s)$.  Then $\tilde{f}(0) = 0$ and $f(0) = n_s$
for some $n_s
> 0$.  Also $f'(0) = \tilde{f}'(0) = c_s$.  Recall that $T^\e_{33}$
is the operator on $B(V_3)$ defined by
   \begin{equation} \label{Te33}
   T^\e_{33}\phi(t_0) = \int \mathbf{1}_{V_3}(t) \tilde{p}^\e(t|t_0)\phi(t)\,dt,
   \quad \quad t_0 \in V_3
   \end{equation}
where
  \begin{equation} \label{tildep}
  \tilde{p}^\e(t|t_0) = \sum_{n \in \mathbb{Z}} p^\e(t+n|t_0)
 \end{equation}
We then extend $T_{33}^\e$ to an operator on $B(\mathbb{R})$ via
  \begin{equation} \label{Te332}
   T^\e_{33}\phi(t_0) = \mathbf{1}_{V_3}(t_0)\int \mathbf{1}_{V_3}(t) \tilde{p}^\e(t|t_0)\phi(t)\,dt,
   \quad \quad t_0 \in \R
   \end{equation}
and look at the re-scaled version $T_s^\e := (U_\e)^{-1} \circ
T_{33}^\e \circ U_\e$ where $U_\e\phi(x) = \phi(x/\e)$, so that
   \begin{eqnarray*} T_s^\e\phi(t_0) &=& \mathbf{1}_{V_3^\e}(t_0)
  \int \mathbf{1}_{V_3^\e}(t) \phi(t) \e \tilde{p}^\e(\e t|\e
   t_0) dt
   \end{eqnarray*}
with $V_3^\e = V_3/\e = (-\delta_s/\e,\delta_s/\e)$.  For ease of
notation, we drop the subscript $s$ from $c_s$ and $\delta_s$ and
the subscript $3$ from $V_3$ and $V_3^\e$.

From (\ref{tildep}) we get
  $$
  \e\tilde{p}^\e(\e t|\e t_0) = \sum_{m \in \mathbb{Z}} \e p^\e(\e t+n|\e t_0).
  $$
Let
  $$p^\e_{main}(t|t_0) = \e p^\e(\e t + n_s |\e t_0)
  $$
denote the main term in this sum.  Theorem \ref{thm-unif},
together with the limiting behavior $m(\e t_0) \to m(0)$ and
$\sigma(\e t_0)\to \sigma(0)$ and $(f(\e t_0)-n_s)/\e \to c t_0$
implies that
   \begin{equation} \label{pconv}
   p^\e_{main}(t|t_0) \to p_\tau(t-ct_0|0) = \frac{1}{\sqrt{2 \pi}
   \sigma_\tau(0)}e^{-(t-ct_0)^2/2\sigma_\tau^2(0)}
   \end{equation}
in some sense.  (For details of this calculation see equation
(\ref{pmainlim}) later.)  This suggests that in some sense,
$T_s^\e \rightarrow T_s$ as $\e \to 0$ where $T_s$ is the operator
with kernel $p_\tau(t-ct_0|0)$.

In order to make this precise, for $k  \in \R$ define $\|\phi\|_k
= \sup\{|\phi(x)|e^{-kx^2}:  x \in \R\}$ and $W_k = \{\phi: \R \to
\R: \phi \mbox{ is measurable and } \|\phi\|_k < \infty\}$.  Then
$W_k$ with the norm $\|\cdot\|_k$ is a Banach space.  Let ${\cal
L}(W_k)$ denote the set of all bounded linear operators $T$ on $W_k$
with operator norm
$$\|T\|_k = \sup\{\|T\phi\|_k: \phi
\in W_k \mbox{ and }\|\phi\|_k \le 1\}.
$$

\begin{proposition} \label{prop-Teps}
For all $k>0$ sufficiently small, we have $T_s^\e = T_s + O(\e)$ in ${\cal L}(W_k)$.
\end{proposition}

\n{\bf Proof.}  Define the operators $T_{main}^\e$, $T^\e_{cut}$
with the following kernels:
   \begin{eqnarray*}
    T_{main}^\e & \leftrightarrow & \mathbf{1}_{V^\e}(t)\mathbf{1}_{V^\e}(t_0) p_{main}^\e (t| t_0)\\[1ex]
   T^\e_{cut} & \leftrightarrow & \mathbf{1}_{V^\e}(t)\mathbf{1}_{V^\e}(t_0) p_\tau(t-ct_0|0)
    \end{eqnarray*}
and write
\begin{eqnarray*}
T_s^\e - T_s &=& (T_s^\e - T_{main}^\e) + (T_{main}^\e-T^\e_{cut}) + (T^\e_{cut}-T_s) \\
&=& I + II + III
\end{eqnarray*}
The proof will consist of bounding $I, II, III$ as operators on $W_k$ for small enough $k$.

\s

To bound $I$, we note that for $\phi \in W_k$ we have
   \begin{align*}
   |(T_s^\e-T_{main}^s)\phi(t_0)|
   & =  \left| \sum_{n \neq n_s}\int \mathbf{1}_V(\e t_0)\mathbf{1}_V(\e t)\phi(t) \e p^\e(\e t+n|\e t_0)\,dt\right|\\
   & \le   \|\phi\|_k e^{k \delta^2/\e^2}\mathbf{1}_V(\e t_0)\left( \sum_{n \neq n_s}\int \mathbf{1}_V(\e t) \e p^\e(\e t+n|\e
   t_0)\,dt\right)\\
      & =  \|\phi\|_k e^{k \delta^2/\e^2}\mathbf{1}_V(\e t_0)P^{\e t_0}\left(\tau^\e \in \bigcup_{n \neq n_s} V+n
   \right).
   \end{align*}
Now $|\e t_0| < \delta$ implies $|f(\e t_0) -n_s| < \delta$, so
that
 $$
 \left\{\tau^\e \in \bigcup_{n \neq n_s} V+n
   \right\} \subset \{|\tau^\e - f(\e t_0)| \ge 1-2 \delta\}.
   $$
We can assume without loss of generality that $\delta$ was chosen
small enough so that $\delta_1:=1-2\delta > 0$.  Then Proposition
\ref{deviations} gives
 \begin{eqnarray*}
   |(T_s^\e-T_{main}^s)\phi(t_0)|
    & \le  &
    \|\phi\|_k e^{k \delta^2/\e^2}P^{\e t_0}\left(|\tau^\e - f(\e t_0)| \ge 1-2 \delta
   \right)\\
   & \le &  \|\phi\|_k e^{k \delta^2/\e^2} M_{\delta_1} \e \, e^{-K_{\delta_1}/\e^2},
   \end{eqnarray*}
and so
  $$
  \|T_s^\e - T_m^\e\|_k \le M_{\delta_1} \e \,  e^{-(K_{\delta_1}-k)/\e^2}.
  $$
This is at most $O(\e)$ as long as $k \le K_{\delta_1}$.

\s

The calculation giving an $O(\e)$ bound for $III$ concerns the
effect of the cutoffs $\mathbf{1}_{V^\e}(t_0)\mathbf{1}_{V^\e}(t)$
on a Gaussian kernel.  This is the same as in the proof of
Equation (17) in \cite{jm} and is omitted.

\s

For $II$ we have
  \begin{eqnarray}
  p_{main}^\e(t|t_0) & =& \e p^\e(f(0)+\e t|\e t_0) \nonumber \\
    & =& \e p^\e(f(\e t_0) + \e t -(f(\e t_0)-f(0))  |\e t_0)  \nonumber \\
     & =& \e p^\e( f(\e t_0) + \e t - \tilde{f}(\e t_0) |\e t_0)  \nonumber \\
       &=& p_\tau(t-ct_0|0) \nonumber \\
       & & + \big[\e p^\e(\e (t - \tilde{f}^\e(t_0)) + f(\e t_0)|\e t_0)-  p_\tau(t-\tilde{f}^\e(t_0)|\e t_0)\big]  \nonumber \\
     & & + \big[p_\tau(t-\tilde{f}^\e(t_0)|\e t_0) - p_\tau(t-ct_0|0)\big]  \nonumber  \\
     &\equiv & p_\tau(t-ct_0|0) + r_1^\e(t|t_0) + r_2^\e(t|t_0) \label{pmainlim}
     \end{eqnarray}
where $\tilde{f}^\e(t_0) = \e^{-1} \tilde{f}(\e t_0)$. Define
operators $T_{r_i}^\e$ with kernels
$$\mathbf{1}_{V^\e}(t)\mathbf{1}_{V^\e}(t_0) r_i^\e(t|t_0),$$
$i=1,2$. The operator $T^\e_{r_2}$ deals with the effect of
changing the mean $\tilde{f}^\e(t_0)$ and standard deviation
$\sigma_\tau(\e t_0)$ of a Gaussian kernel, and Equation (14) in
\cite{jm} gives an $O(\e)$ bound for $\|T_{r_2}^\e\|_k$ for all
sufficiently small $k>0$.   For $T^\e_{r_1}$, suppose $k > 0$ and
$\phi \in W_k$.  We have
 \begin{align*}
  |T^\e_{r_1} \phi(t_0)| & = \left|\int \mathbf{1}_{V^\e}(t)\mathbf{1}_{V^\e}(t_0) r_1^\e(t|t_0)
  \phi(t)\,dt\right|\\
   & =   \left|\int \mathbf{1}_{V^\e}(t)\mathbf{1}_{V^\e}(t_0) \right. \\
    &  \quad \times \left.[\e p^\e(\e (t - \tilde{f}^\e(t_0)) + f(\e t_0)|\e t_0) -  p_\tau(t-\tilde{f}^\e(t_0)|\e t_0)]
   \phi(t)\,dt\right|.
     \end{align*}
In order to apply Theorem 1 we need to restrict to $|\e
t-\tilde{f}(\e t_0)|$ sufficiently small.  This can be achieved
for $t_0,t \in V^\e$ by choosing $\delta_s$ sufficiently small in
Proposition \ref{regions}. Then
  \begin{eqnarray*}
  |T^\e_{r_1} \phi(t_0)|
   & \le &  K \e \int \mathbf{1}_{V^\e}(t)\mathbf{1}_{V^\e}(t_0)
    e^{-(t - \tilde{f}^\e(t_0))^2/2\sigma_1^2}
    |\phi(t)|\,dt \\
        & \le & K \e \|\phi\|_k \mathbf{1}_{V^\e}(t_0)\int
    e^{-(t - \tilde{f}^\e(t_0))^2/2\sigma_1^2} e^{kt^2}\,dt\\
    & = & K \e \|\phi\|_k \mathbf{1}_{V^\e}(t_0) \frac{\sqrt{2 \pi} \sigma_1}{\sqrt{1-2k\sigma_1^2} }\exp\left\{\frac{k
    (\tilde{f}^\e(t_0))^2}{1-2k\sigma_1^2}\right\}.
  \end{eqnarray*}
Furthermore by shrinking $V$ if necessary, we can find $c_1 < 1$
such that $|\tilde{f}^\e(t_0)| \le c_1 |t_0|$ for $t_0 \in V^\e$,
and then
   $$
  |T^\e_{r_1} \phi(t_0)|
   \le  K \e \|\phi\|_k \mathbf{1}_{V^\e}(t_0) \frac{\sqrt{2 \pi} \sigma_1}{\sqrt{1-2k\sigma_1^2} }\exp\left\{\frac{k
    c_1^2 t_0^2}{1-2k\sigma_1^2}\right\}.
     $$
It follows that
     $$
  \|T^\e_{r_1} \|_k
   \le  K \e \frac{\sqrt{2 \pi} \sigma_1}{\sqrt{1-2k\sigma_1^2} }
     $$ so long as $c_1^2 < 1-2k
\sigma_1^2$, that is, $k < (1-c_1^2)/2\sigma_1^2$.  This completes
the proof of Proposition \ref{prop-Teps}.  \hfill $\Box$

\s

Now $T_s$ is the transition operator for the Markov chain $X_n =
cX_{n-1} + \chi_n$, and has eigenvalues $c^n$, $n \ge 0$.  It
follows from Proposition \ref{prop-Teps} together with standard
perturbation results for linear operators (see Kato \cite{tk})
that the eigenvalues of $T_s^\e$ acting on $W_k$ are of the form
$c^n +O(\e)$.  Finally, since the operator $T_s^\e$ involves the
indicator functions $\mathbf{1}_{V^\e}(t_0)\mathbf{1}_{V^\e}(t)$
the eigenvalues and eigenfunctions do not depend on the value of
$k$, and in particular the non-zero eigenvalues of $T_s^\e$ acting
on $W_k$ coincide with the non-zero eigenvalues of $T_{33}^\e$
acting on $B(V)$.  This completes the proof of Proposition
\ref{prop-stabev}.

\subsection{Behavior near an unstable fixed point}  \label{sec-unstable}

Here $T_{1}^\e$ is the restriction of the transition operator
$T^\e$ to a neighborhood $V_1$ of the unstable fixed point
$\theta_u$, and $f'(\theta_u) = c_u$ with $|c_u| > 1$.  The main
result of this section is

\begin{proposition} \label{prop-unstabev}  Every non-zero eigenvalue
of $T_{11}^\e$ is of the form
$$|c_u|^{-1}c_u^{-n} + O(\e)$$
as $\e \to 0$ for some $n \ge 0$.
\end{proposition}

Here we localize in the neighborhood $V_1$ of the unstable fixed
point and replace $T_s^\e$ and $T_s$ in Section \ref{sec-local}
with $T_u^\e$ and $T_u$ defined in
the same way, but using $V_1$ in place of $V_3$. We obtain a perturbation result similar to
Proposition \ref{prop-Teps}, but in a class of spaces of functions
with exponential decay. Note that all calculations in the previous
section before the proof of Proposition \ref{prop-Teps} are valid
if $|c| > 1$ as well.

\begin{proposition}
For all $k>0$ sufficiently small, we have $T_u^\e = T_u + O(\e)$ in ${\cal L}(W_{-k})$.
\end{proposition}

\n{\bf Proof.} Define the operators $T_{main}^\e$ and $T_{cut}^\e$
as in the proof of Proposition \ref{prop-Teps} and decompose
$T_u^\e - T_u$ in the same way as $ I+II+III$. The bound on $I$ is
obtained in essentially the same as in Section \ref{sec-local};
the only difference is that now $\e t_0 \in V_1$ implies $|f(\e
t_0)-n_u| \le \hat{\delta}$ for some $\hat{\delta}
> 0$, so that we need $\delta_1:= 1-\delta - \hat{\delta} > 0$.
The bounds on $III$ and the bound on $T_{r_2}^\e$ in the
decomposition of $II$ involve Gaussian kernels and are the same as
those used in the proof of Theorem 9 in \cite{jm}. Finally, by
shrinking $V_1$ if necessary, we can find $c_1 > 1$ such that
$|f^\e(t_0)| > c_1 t_0$ for all $\e t_0 \in V_1$.  By a similar
argument to that given in Section \ref{sec-local} above, for $k
> 0$ and $\phi \in W_{-k}$ we obtain
    $$
  |T^\e_{r_1} \phi(t_0)|
   \le  K \e \|\phi\|_{-k} \mathbf{1}_{V^\e}(t_0) \frac{\sqrt{2 \pi} \sigma_1}{\sqrt{1+2k\sigma_1^2} }\exp\left\{\frac{-k
    c_1^2 t_0^2}{1+2k\sigma_1^2}\right\}.
     $$
It follows that
   $$
  \|T^\e_{r_1} \|_{-k}
   \le  K \e \frac{\sqrt{2 \pi} \sigma_1}{\sqrt{1+2k\sigma_1^2} }
     $$
so long as $c_1^2 \ge 1+2 k \sigma_1^2$, that is, $k \le
(c_1^2-1)/2\sigma_1^2$.  \hfill $\Box$

\s

Since $T_u$ has eigenvalues $|c|^{-1}c^{-n}$ if $|c| > 1$ (see
Section 5 of \cite{jm}), the proof of Proposition
\ref{prop-unstabev} again follows from standard perturbation
arguments.

\section{Proofs for Section \ref{sec-disc}} \label{sec-discproofs}

\n {\bf Proof of Proposition \ref{deviations2}.}  Since $D$ is finite it suffices to prove the
result separately for each $t_0 \in D$.  Suppose first that $t_0
\in D$ satisfies {\bf (C')(i)}, so that $f^*(t_0) = f(t_0)$. Since
$\xi$ is a continuous function of both arguments and $g$ is
continuous at $s=t_i$, there exist $\tilde{\delta}, \delta_1 > 0$
and $\delta_2 \in(0,\delta]$ such that
   $$
     \xi(s|t) < g(s) - \delta_1 \quad \mbox{ for } t \le s \le f(t_0) - \delta
     $$
and
     $$
     \xi(s|t) > g(s) + \delta_1 \quad \mbox{ for } s  =
     f(t_0)+\delta_2.
     $$
whenever $|t-t_0| < \tilde{\delta}$.  As in the proof of
Proposition \ref{deviations}, these two inequalities imply that
  \begin{eqnarray*}
  \PP^{t}(|\tau^\e - f(t_0)| > \delta)
    & \le & \PP^{t}\left(\sup_{t \le s \le f(t_0)+\delta_2}|X^\e(s) - \xi(s|t)| \ge \delta_1\right).
  \end{eqnarray*}
If instead $t_0 \in D$ satisfies {\bf (C')(ii)}, then again using
the continuity of $\xi$ and $g$, there exist $\tilde{\delta} > 0,
\delta_1 > 0$ and $\delta_2 \in(0,\delta]$ such that
   $$
     \xi(s|t) < g(s) - \delta_1 \quad \mbox{ for } s \in[t,f(t_0) -
     \delta) \cup (f(t_0)+\delta,f^*(t_0) - \delta)
     $$
and
     $$
     \xi(s|t) > g(s) + \delta_1 \quad \mbox{ for } s  =
     f^*(t_0)+\delta_2.
     $$
whenever $|t-t_0| < \tilde{\delta}$.  These two inequalities imply
that
  \begin{eqnarray*}
  \PP^{t}(d(\tau^\e,\{f(t_0),f^*(t_0)\}) > \delta)
    & \le & \PP^{t}\big(\{\tau^\e < f(t_0)-\delta\} \\
    & & \hspace{7ex} \cup \{f(t_0)+\delta < \tau^\e < f^*(t_0)-\delta\}\\
      & & \hspace{7ex} \cup\, \{\tau^\e > f^*(t_0)+\delta_2\}\big)\\
    & \le &  \PP^{t}\left(\sup_{t \le s \le f^*(t_0)+\delta_2}
     |X^\e(s) - \xi(s|t)| \ge \delta_1\right).
     \end{eqnarray*}
In either case the proof is completed using the same method as in
the proof of Proposition \ref{deviations}.  \hfill $\Box$

\s

The following Lemma will be used several times in the proof of
Theorem \ref{thm-discspec}.  Recall that $\tilde{f} = f$ mod 1.

\begin{lemma} \label{lem-1} Assume that $f$ satisfies {\bf (A)}, {\bf (B')} and {\bf (C')}.
For $\theta \in S$ suppose $\tilde{f}^i(\theta) \not \in D$ for $0
\le i \le n-1$ and let $V$ be any neighborhood of
$\tilde{f}^n(\theta)$. There is a neighborhood $U$ of $\theta$ and
constants $K$ and $M$ such that
   $$\PP^\psi(\Theta^\e_n \not\in V) \le M\e e^{-K/\e^2}
   $$
for $\psi \in U$.
  \end{lemma}

\n{\bf Proof.}  For $n = 1$ this is a simple consequence of
Proposition 2, and the proof for general $n$ follows by a simple
inductive argument.  \hfill $\Box$

\s

\n{\bf Proof of Theorem \ref{thm-discspec}.}  By {\bf (D1)}, for
any $\delta_1 > 0$ we can choose an open set $V_2$ with $P \subset
V_2 \subset \bigcup_{i=1}^k \overline{B}(\theta_1,\delta_1)
\subset \mathbb{S} \backslash D$ so that
$d(\tilde{f}(\theta),V_2^c) > \delta$ for all $\theta \in V_2$ and
some constant $\delta >0$.   Let $V_1 = \mathbb{S}\backslash V_2$.
Then, similarly to the proof of Theorem \ref{thm-outrans}, we can
use the decomposition $S = V_1 \cup V_2$ to write
    $$
T^\e = \left[ \begin{array}{cc} T_{11}^\e & T_{12}^\e \\
      T_{21}^\e & T_{22}^\e \end{array} \right].
$$
By Proposition 2 we have
   $$
   \sup_{\theta \in V_2} \PP^\theta(\Theta^\e_1 \not\in V_2) \le M_1
\e e^{-K_1/\e^2}
   $$
for some $K_1$ and $M_1$, and therefore $\|T^\e_{21}\|_\infty \le
M_1 \e e^{-K_1/\e^2}$. Moreover, for any $n \ge 1$ we have
    \begin{equation} \label{inV2}
   \sup_{\theta \in V_2} \PP^\theta(\Theta^\e_n \not\in V_2) \le n M_1
\e e^{-K_1/\e^2}.
   \end{equation}
Next we estimate $\|(T^\e_{11})^N\|_\infty$ for sufficiently large
$N$.  Using {\bf (D3)}, we can write the compact set $V_1$ as the
disjoint union $A_0 \cup A_1 \cup \cdots \cup A_{\ell-1} \cup B$
where $A_i = \tilde{f}^{-i}(D)$.  (The disjointness of the $A_i$
follows from the fact that $f^j(\theta) \not\in D$ for $\theta \in
D$ and $j \ge 1$.)  Notice that for $\theta \in B$ we have
$f^i(\theta) \not\in D$ for all $i \ge 0$, and that $E \subset B$.

For each $\theta \in B$, the condition {\bf (D2)} implies there
exists $n_\theta$ such that $f^{n_\theta}(\theta) \in V_2$ and
$f^i(\theta) \not \in D$ for $0 \le i \le n_\theta - 1$.  By Lemma
\ref{lem-1} with $V = V_2$ there exist a neighborhood $U_\theta$
of $\theta$ and constants $K_\theta$ and $M_\theta$ such that
\begin{equation} \label{fromB}
   \sup_{\psi \in U_\theta} \PP^\psi(\Theta^\e_{n_\theta} \not\in V_2) \le M_\theta
\e e^{-K_\theta/\e^2}.
  \end{equation}
Now suppose $\theta \in A_0 = D$.  Then $f(\theta)$ and
$f^*(\theta)$ are both in $B$.  By Proposition \ref{deviations2}
there exist a neighborhood $U_\theta$ of $\theta$ and constants
$\widetilde{K}_\theta$ and $\widetilde{M}_\theta$ such that
   \begin{equation} \label{DtoB}
   \sup_{\psi \in U_\theta} \PP^\psi(\Theta^\e_1 \not\in (U_{f(\theta)} \cup U_{f^*(\theta)})) \le \widetilde{M}_\theta
\e e^{-\widetilde{K}_\theta/\e^2}.
   \end{equation}
Define $n_\theta = \max(n_{f(\theta)},n_{f^*(\theta)})+1$.
Combining (\ref{DtoB}) and (\ref{fromB}) and (\ref{inV2}) gives
   \begin{equation} \label{DtoV2}
   \sup_{\psi \in U_\theta} \PP^\psi(\Theta^\e_{n_\theta} \not\in V_2) \le M_\theta
\e e^{-K_\theta/\e^2}
   \end{equation}
for some $K_\theta$ and $M_\theta$.  Finally suppose that $\theta
\in A_i$ for $i \ge 1$.  Then $f^j(\theta) \not\in D$ for $0 \le j
< i$ and $f^i(\theta) \in D$. By Lemma \ref{lem-1} with $V =
U_{f^i(\theta)}$ there exist a neighborhood $U_\theta$ of $\theta$
and constants $\widetilde{K}_\theta$ and $\widetilde{M}_\theta$
such that
    \begin{equation} \label{AitoD}
   \sup_{t \in U_\theta} \PP^t(\Theta^\e_i \not\in U_{f^i(\theta)}) \le \widetilde{M}_\theta
\e e^{-\widetilde{K}_\theta/\e^2}.
   \end{equation}
Define $n_\theta = i+ n_{f^i(\theta)}$.  Combining (\ref{AitoD})
and (\ref{DtoV2}) gives
     \begin{equation} \label{AitoV2}
   \sup_{\psi \in U_\theta} \PP^\psi(\Theta^\e_{n_\theta} \not\in V_2) \le M_\theta
\e e^{-K_\theta/\e^2}.
   \end{equation}
for some $K_\theta$ and $M_\theta$.

We have constructed an open cover $\{U_\theta: \theta \in B \cup
A_0 \cup \cdots \cup A_{\ell - 1}\}$ of the compact set $V_1$.
Passing to a finite subcover $\{U_{\theta_j}\}$ and letting $N =
\max_j n_{\theta_j}$ and using (\ref{inV2}) with $n =
N-n_{\theta_j}$ together with (\ref{fromB}) and (\ref{DtoV2}) and
(\ref{AitoV2}) gives
    \begin{equation} \label{V1toV2}
   \sup_{t \in V_1} \PP^t(\Theta^\e_N \not\in V_2) \le M
\e e^{-K/\e^2}
   \end{equation}
for some $K$ and $M$, and thus $\|(T^\e_{11})^N\|_\infty \le M \e
e^{-K/\e^2}$.

It remains only to describe the eigenvalues of $T^\e_{22}$ and
this can be done using exactly the same methods used for Theorem
\ref{thm-outrans} and Remark \ref{rmk-periodp}.  \hfill $\Box$

\end{document}